\documentclass[12pt,a4paper,draft]{article}

\usepackage[reqno]{amsmath}
\usepackage{amssymb,amsthm,upref,amscd}
\usepackage[T1]{fontenc}
\usepackage{times}

\newtheorem{theorem}{Theorem}[section]
\newtheorem{corollary}[theorem]{Corollary}

\newtheorem{example}[theorem]{Example}

\newtheorem{lemma}[theorem]{Lemma}
\newtheorem{proposition}[theorem]{Proposition}

\theoremstyle{definition} \theoremstyle{remark}
\numberwithin{equation}{section}

\begin{document}

\title{The existence of solutions of Schr\"odinger equations with essence resonance
\\
[3mm] {\footnotesize Dedicated to Professor Shujie Li on the occasion of his
85th birthday }}
\author{Chong Li$^{\mathrm{a,b,c,d}}$\thanks{%
The author is supported by NSFC(11871066), AMSS(E3550105). Email:
lichong@amss.ac.cn.} \\
\\
{\small $^{\mathrm{a}}$State Key Laboratory of Mathematical Sciences,
Academy of Mathematics and Systems Science, }\\
{\small Chinese Academy of Sciences, Beijing 100190, Peoples Republic of
China}\\
{\small $^{\mathrm{b}}$Institute of Mathematics, Academy of Mathematics and
Systems Science, Chinese Academy of Sciences}\\
{\small $^{\mathrm{c}}$School of Mathematical Sciences, University of
Chinese Academy of Sciences, Beijing 100049, China}\\
{\small $^{\mathrm{d}}$Center for Excellence in Mathematical Sciences,
Chinese Academy of Sciences }}
\date{}
\maketitle

\begin{abstract}
The current paper investigates a class of asymptotically linear
Schr\"odinger equations with the nonlinear term $f\in C\left( \Bbb{R},\Bbb{R}%
\right) $ s.t. $\lim\limits_{\left| t\right| \rightarrow \infty }\frac{%
f\left( t\right) }t=\sigma _0$, where $\sigma _0$ stands for the threshold
of essential spectrum of Schr\"odinger operator. Clearly the Palais-Smale
condition fails to hold in this case. Especially under the hypothesis $%
\left( V_2\right)$, the lack of compactness occurs at the interaction
between nonlinear term and continuum spectrum. For this reason, we introduce
a bootstrap iteration approach for elliptic equation on $\Bbb{R}^N$. The
iteration is self-contained and can be regarded as a generalization of
Agmon-Douglis-Nirenberg theorem(see \cite{ADN1},\cite{ADN2}). The proof
characterizes iteration steps independent of the choice of the parameter $%
\lambda $, which are indeed manipulated by intrinsic natures of potentials
and nonlinear terms, and furthermore presents precise estimates for
asymptotically linear functions or continuous nonlinear terms restricted on
a bounded domain $\Omega $ with smooth boundary $\partial \Omega $ in $\Bbb{R%
}^N$. Additionally, a comparison theorem for the spectrum of Schr\"odinger
operator is also established in this paper. With above preparations, we can
get a nontrivial solution without mountain pass geometry, and more
importantly make an explicit description of nondegeneracy of solutions with
monotonicity hypothesis. \bigskip


\end{abstract}


\medskip

\section{Introduction}

This paper is mainly concerned with nonlinear Schr\"odinger equation:
\begin{equation}
\left\{
\begin{array}{cc}
-\Delta u+V\left( x\right) u=\lambda u+g\left( u\right) \text{\textrm{,}} &
x\in \Bbb{R}^N\text{\textrm{,}} \\
u\left( x\right) \rightarrow 0\text{\textrm{,}} & \text{\textrm{as }}\left|
x\right| \rightarrow \infty \text{\textrm{,}}
\end{array}
\right.  \label{eq1}
\end{equation}
arising from study of standing wave solutions of time-dependent nonlinear
Schr\"odinger equations, where $G\left( u\right) =\int_0^ug\left( s\right)
ds $, and $g\left( s\right) $ is a Carath\'eodory function on $\Bbb{R}$ .
Let $V $ be a real potential and the form domain of $-\Delta +V$ is the
Hilbert space $H^1\left( \Bbb{R}^N\right) $ equipped with the norm
\[
\left\| u\right\| _{H^1}=\left( \int_{\Bbb{R}^N}\left| \nabla u\right|
^2dx+\int_{\Bbb{R}^N}\left| u\right| ^2dx\right) ^{\frac 12}
\]
and the inner product
\[
\left\langle u,v\right\rangle =\int_{\Bbb{R}^N}\left( \nabla u\nabla
v+uv\right) dx\text{\textrm{.}}
\]
The corresponding energetic functional of $\left( \ref{eq1}\right) $ is of
the form as:
\[
J_\lambda \left( u\right) =\frac 12\int_{\Bbb{R}^N}\left( \left| \nabla
u\right| ^2+V\left( x\right) u^2\right) dx-\frac \lambda 2\int_{\Bbb{R}%
^N}u^2dx-\int_{\Bbb{R}^N}G\left( u\right) dx\text{\textrm{, }}u\in H^1\left(
\Bbb{R}^N\right) \text{\textrm{,}}
\]
and from the variational point of view, the solutions of $\left( \ref{eq1}%
\right) $ are the critical points of $J_\lambda $.

For the sake of convenience, throughout the paper we always assume that the
linear potential $V$ is a Kato-Rellich potential(A real potential $V\left(
x\right) $ is called Kato-Rellich(K-R)potential if we can decompose $%
V=V_1+V_2$, $V_1\in L^p\left( \Bbb{R}^N\right) $, $V_2$ $\in L^\infty \left(
\Bbb{R}^N\right) $, with $p=2$ if $N\leq 3$, $p>2$ if $N=4$ and $p>\frac N2$
if $N\geq 5$, see \cite{HS}). Let $A=-\Delta +V$ and assume:

$\left( V_1\right) $ $\sigma _{\text{\textrm{dis}}}\left( A\right) \cap
\left[ \mu _1,\sigma _0\right) =\left\{ \mu _j\right\} _{j=1}^k$, $k\geq 3$,
s.t. $\mu _0:=-\infty <\mu _1:=\inf \sigma \left( A\right) <\mu _2\leq
\cdots \leq \mu _j\leq \cdots \leq \mu _k<\sigma _0:=\inf \sigma _{\mathrm{%
ess}}\left( A\right) $;

$\left( V_2\right) $ $\lim\limits_{\left| x\right| \rightarrow \infty
}V\left( x\right) =\sigma _0$.

Let $g\in C^1\left( \Bbb{R}^N\times \Bbb{R},\Bbb{R}\right) $, $g\left(
0\right) =0$, and make the following hypotheses:

$\left( g_1\right) $ $g_\infty :=\lim\limits_{\left| t\right| \rightarrow
\infty }\frac{g\left( x,t\right) }t=0$ uniformly with respect to $x\in \Bbb{R%
}^N$;

$\left( g_2\right) $ $g_0:=\lim\limits_{t\rightarrow 0}\frac{g\left(
x,t\right) }t$ uniformly with respect to $x\in \Bbb{R}^N$.

$\left( g_3\right) $ $\frac{g\left( x,t\right) }t\leq 0$ for $\forall s\in
\Bbb{R\backslash }\left\{ 0\right\} $, $\forall t\neq 0$.

For the case $V\left( x\right) =$ constant $>0$ and the nonlinear term $%
f\left( x,t\right) $ is $1$-periodic, Jeanjean\cite{J} obtained a positive
solution by a generic theorem, and monotonicity hypothesis was employed by
\cite{LZ} to make an argument on the existence of a positive solution. By
supposing $V\left( x\right) \geq \alpha >0$ and applying mountain pass
theorem, Jeanjean and Tanaka \cite{JT} got a positive solution via Ekeland's
principle. For a peculiar class of potentials $V\left( x\right) $, s.t.,
there exists some $R>0$ such that
\[
V\left( x\right) \leq \sigma _0,\text{ a.e. }x\in B_R:=\left\{ x\in \Bbb{R}%
^N:\left| x\right| <R\right\} ,
\]
and
\[
\int_{B_R}\left( V\left( x\right) -\sigma _0\right) +\int_{B_R^c\cap \left\{
V>\sigma _0\right\} }\left( V\left( x\right) -\sigma _0\right) \leq 0,
\]
Song\cite{S} derived a positive solution. It should be noted that aiming to
seek a nontrivial solution mountain-pass geometry of energetic functional is
assumed by \cite{CT},\cite{J},\cite{JT},\cite{LZ},\cite{S},\cite{SZ} without
exception. Furthermore, with the hypothesis of radical symmetry on $V\left(
x\right) $ and $f\left( x,t\right) $, \cite{MS} and \cite{W} explored
respectively the indefinite and singular potential cases.

A class of asymptotically linear problems are called strong resonance, in
which Landesman-Lazer (L-L) condition is dropped so that the Palais-Smale
condition does not hold.

Evidently essence resonance(we say essence resonance if the nonlinear term $%
f\in C\left( \Bbb{R}^N\times \Bbb{R},\Bbb{R}\right) $ s.t. $%
\lim\limits_{\left| t\right| \rightarrow \infty }\frac{f\left( x,t\right) }%
t=\sigma _0$ uniformly with respect to $x\in \Bbb{R}^N$, where $\sigma _0$
denotes the infimum of essential spectrum of Schr\"odinger operator) can be
treated as a kind of sophisticated strong resonance. Based on the absence of
compactness for essence resonance problems, \textbf{a natural question which
arises is: Can we still derive the existence of solutions of }$\left( \ref
{eq1}\right) $\textbf{\ }for $\lambda \in \sigma _{\text{\textrm{ess}}%
}\left( A\right) $ \textbf{if the hypothesis that the energetic functional
possess the mountain-pass geometry is eliminated?}

The current paper intends to approach the subject by focusing on nonlinear
Schr\"odinger equations interacting with essential spectrum, and to do that
we need some additional assumptions:

$\left( g_4\right) $ $g_t^{\prime }\left( x,t\right) >\frac{g\left(
x,t\right) }t$, $\forall t\in \Bbb{R\backslash }\left\{ 0\right\} $, $%
\forall x\in \Bbb{R}^N$.

$\left( g_5\right) $ There exists a $\beta >0$, s.t., $\forall s_1$,$s_2\in
\Bbb{R}$, $\forall x\in \Bbb{R}^N$,
\[
\left| g\left( x,s_1\right) -g\left( x,s_2\right) \right| \leq \beta \left|
s_1-s_2\right| ,
\]

$\left( g_6\right) $ $\lambda +g_0\leq \mu _i$ for $i\leq k$.

$\left( g_7\right) $ $\lim\limits_{t\rightarrow \infty }\frac{G\left(
x,t\right) }{t^2}=0$ uniformly with respect to $x\in \Bbb{R}^N$.

$\left( g_8\right) $ $g_t^{\prime }\left( x,t\right) \leq 0$ for $t\neq 0$.

Stemming from the compactness argument by \cite{LiLi}, we proceed to
validate the feasibility of Palais-Smale condition for the nonresonance case
under a more general assumption $\left( g_5\right) $ in Section 3. The
outline of the proof of the existence assertion offered by Section 5 can be
summarized as follows:

$\left( i\right) $ Construct a linking for the nonresonance case $\lambda
\notin \sigma \left( A\right) $.

$\left( ii\right) $ Verify that for $a\in \Bbb{R}$, $a<\sigma _0$, $\left[
a,\sigma _0\right) \cap \sigma \left( A\right) =\varnothing $, the critical
value $c_\lambda $ of the functional $J_\lambda $ yielded by linking
approach is bounded for $\lambda \in \left[ a,\sigma _0\right) $, and
furthermore $\left\| u_\lambda \right\| _{H^1}$ does as well for the
corresponding critical point $u_\lambda $ with $c_\lambda =J_\lambda \left(
u_\lambda \right) >0$.

$\left( iii\right) $ Demonstrate that $\left\| u_\lambda \right\| _{L^\infty
}$ is indeed bounded by invoking bootstrap iteration for Schr\"odinger
equations built up in Section 2.

$\left( iv\right) $ Based on the above deductions and by confirming a
comparison theorem for the spectrum of Schr\"odinger operator in Section 4,
we elucidate that for an arbitrarily picked sequence $\lambda _n\rightarrow
\sigma _0$, $\left\{ u_{\lambda _n}\right\} _{n=1}^\infty $ converges to a
nonzero critical point $u_{\sigma _0}$ of $J_{\sigma _0}$, and besides,
adding some appropriate hypotheses we can show the nondegeneracy of the
solution $u_{\sigma _0}$.

Compared to the treatment on the class of superquadratic functions, we need
one more trick to address the problems on establishing a topological link
for asymptotically linear equations, and this point will be embodied in the
current paper(the hypothesis $\left( g_7\right) $ is easy to be validated
and also essential during the elaboration process). In our opinions, the
significance of the geometry construction is expressed as follows: we
initially obtain the boundedness of critical values produced by linking for $%
J_\lambda $ with $\lambda \in \left[ a,\sigma _0\right) $, $\left[ a,\sigma
_0\right) \cap \sigma \left( A\right) =\varnothing $, and then employing a
compactness argument furnished by \cite{LiLi} coupling with
concentration-compactness argument, we can derive that the solution set
composed of the critical points yielded by linking approach is bounded in $%
H^1\left( \Bbb{R}^N\right) $. Denote by $\bigcup\limits_{\lambda \in \left[
a,\sigma _0\right) }\left\{ u_\lambda \right\} $ the solution set(indeed we
can further get $0\notin \bigcup\limits_{\lambda \in \left[ a,\sigma
_0\right) }\left\{ u_\lambda \right\} $ under the hypotheses $\left(
g_4\right) \left( g_6\right) $). Given a bounded open domain $\Omega \subset
\Bbb{R}^N$ with smooth boundary $\partial \Omega $ and $\frac Nm>q\geq 1$,
we call $W^{m,q}\left( \Omega \right) \hookrightarrow L^{q^{*}}\left( \Omega
\right) $ critical embedding, $q^{*}=\frac{qN}{N-mq}$. It is well-known that
there exists a constant $C=C\left( m,q,N\right) >0$, $\left\| u\right\|
_{W^{m,q}\left( \Omega \right) }\geq C\left\| u\right\| _{L^{q^{*}}\left(
\Omega \right) }$ for arbitrary $u\in W^{m,q}\left( \Omega \right) $, and
this shows that $C$ is independent of the choice of $\Omega $, and moreover
the same is true for $\Omega =\Bbb{R}^N$(see \cite{AF},\cite{B}). In this
paper we adopt an alternative bootstrap iteration(Under the hypotheses of
Theorem\ref{thm2.2}, the iteration consists of critical embeddings for the
case $V_1=0$ and otherwise even more intricate(see Section 2 for details)).
Unlike the usual way, the iterative approach is especially favorable for
precisely estimating iteration steps, and ultimately gives a uniform
estimation of $\bigcup\limits_{\lambda \in \left[ a,\sigma _0\right)
}\left\{ u_\lambda \right\} $ in $L^\infty \left( \Bbb{R}^N\right) $. We
also show in Section 5 that the assumption $\left( g_6\right) $ is
indispensable for investigating the existence of nonzero solution of $\left(
\ref{eq1}\right) $.

Let $E$ be a Hilbert space and $\Phi \in C^2\left( E,\Bbb{R}\right) $. We
write $m\left( \Phi ,u\right) $ as the sum of the dimensions of the negative
eigenspace of the Hessian of energetic functional at $u\in E$, and also
denote by $M\left( \Phi ,u\right) $ the generalized Morse index of at $u$ to
be the sum of the dimensions of the negative and null eigenspace of the
operator $\Phi ^{\prime \prime }\left( u\right) $ acting on $E$ if $u$ is a
critical point of $\Phi $ on $E$. Hence, our main consequence concerning the
existence of solutions for a class of nonlinear Schr\"odinger equations $%
\left( \ref{eq1}\right) $ with essence resonance reads:

\begin{theorem}
\label{thm1.1} Let $V$ be a K-R potential and suppose $N\neq 2$ for the case
$V_1\neq 0$. Given $\lambda =\sigma _0$, under the hypotheses $\left(
V_1\right) \left( V_2\right) \left( g_1\right) $-$\left( g_7\right) $, if $%
\sigma _0+g_0\leq \mu _i<\mu _k$ or $\sigma _0+g_0<\mu _i$, then $\left( \ref
{eq1}\right) $ possesses a nontrivial solution. In addition, if $\sigma
_0+g_0>\mu _{k-1}$ and $\left( g_8\right) $ holds, then there exists $\delta
_0>0$, for $\lambda \in \left[ \sigma _0-\delta _0,\sigma _0+\delta
_0\right] $, $\left( \ref{eq1}\right) $ admits a nondegenerate solution $%
u_\lambda $ with $m\left( J_\lambda ,u_\lambda \right) =k$.
\end{theorem}

It might just as well assume $m+\mu _1>0$ and $\inf \sigma \left( \widehat{A}%
_m\right) >1$ with $\widehat{A}_m=A+V+2m$ by choosing $m>0$ suitably large.
Clearly, $A_m:=A+m$ is a positive definite self-adjoint operator with $%
D\left( A_m\right) =H^2\left( \Bbb{R}^N\right) $. Denote by
\[
H_m^1\left( \Bbb{R}^N\right) =\left\{ u\in H^1\left( \Bbb{R}^N\right) :\int_{%
\Bbb{R}^N}\left| \nabla u\right| ^2+\left( V+m\right) u^2<\infty \right\}
\]
the Hilbert space equipped with the norm
\[
\left\| u\right\| _m=\left[ \int_{\Bbb{R}^N}\left| \nabla u\right| ^2+\left(
V+m\right) u^2\right] ^{\frac 12}
\]
and the inner product
\[
\left\langle u,v\right\rangle _m=\int_{\Bbb{R}^N}\nabla u\nabla v+\left(
V+m\right) uv.
\]

It is easy to check that $H_m^1\left( \Bbb{R}^N\right) =$ $H^1\left( \Bbb{R}%
^N\right) $ for a K-R potential $V$. Indeed, notice that
\begin{eqnarray}
\frac{\left\| u\right\| _{H^1}^2}2 &\leq &\frac 12\int_{\Bbb{R}^N}\left|
\nabla u\right| ^2+\frac 12\int_{\Bbb{R}^N}\left| \nabla u\right| ^2+2\left(
V+m\right) u^2  \nonumber \\
\ &=&\left\| u\right\| _m^2\leq \left( C\left\| V_1\right\| _{L^p}+\left\|
V_2\right\| _{L^\infty }+m+1\right) \left\| u\right\| _{H^1}^2,  \label{eq2}
\end{eqnarray}
where $p=2$ if $N\leq 3$, and $p>\frac N2$ if $N\geq 4$, so the assertion
follows.

Henceforth, unless otherwise specified, we shall take $X=H_m^1\left( \Bbb{R}%
^N\right) $. On account of the hypothesis, $\sigma \left( A_m\right) \in
\left( 0,+\infty \right) $, alluding to the fact that $A_m^{-1}$ exists. Set
\begin{eqnarray*}
g_{\lambda ,m}\left( u\right) &:&=g_\lambda \left( u\right) +mu\text{\textrm{%
;}} \\
g_\lambda \left( u\right) &:&=\lambda u+g\left( u\right) \text{\textrm{.}}
\end{eqnarray*}
Obviously, in view of the hypotheses, $A_m^{-1}g_{\lambda ,m}\left( u\right)
\in H^2$ for $\forall u\in H^1$. Hence, for fixed $\varphi \in H^1$,
\begin{eqnarray}
\int_{\Bbb{R}^N}g_{\lambda ,m}\left( u\right) \varphi &=&\int_{\Bbb{R}%
^N}A_mA_m^{-1}g_{\lambda ,m}\left( u\right) \cdot \varphi  \nonumber \\
\ &=&\int_{\Bbb{R}^N}\nabla A_m^{-1}g_{\lambda ,m}\left( u\right) \cdot
\nabla \varphi +\left( V+m\right) A_m^{-1}g_{\lambda ,m}\left( u\right)
\cdot \varphi .  \label{eq3}
\end{eqnarray}
Therefore,
\begin{eqnarray}
\left\langle J_\lambda ^{\prime }\left( u\right) ,\varphi \right\rangle _m
&=&\left\langle u-A_m^{-1}g_{\lambda ,m}\left( u\right) ,\varphi
\right\rangle _m  \nonumber \\
\ &=&\int_{\Bbb{R}^N}\nabla u\nabla \varphi +V\left( x\right) u\varphi
-\int_{\Bbb{R}^N}g_\lambda \left( u\right) \varphi ,  \label{eq4}
\end{eqnarray}
yielding $J_\lambda ^{\prime }\left( u\right) =u-A_m^{-1}g_{\lambda
,m}\left( u\right) $.

\section{The global regularity of solutions of Schr\"odinger equations with
Kato-Rellich potentials and a generalization of Agmon-Douglis-Nirenberg
theorem}

This section is devoted to establishing the global regularity of solutions
for the Schr\"odinger equation $\left( \ref{eq1}\right) $ with Kato-Rellich
potentials, and to presenting a corresponding extension of the classical
Agmon-Douglis-Nirenberg theorem. The main result, Theorem\ref{thm2.2}, is
proven via a self-contained bootstrap iteration argument, which constitutes
a principal methodological contribution. The proof's logical flow is
structured as follows. First, the standard Agmon-Douglis-Nirenberg estimates
for elliptic equations on bounded domains are recalled (Proposition\ref
{prop2.1}). This naturally motivates the investigation of analogous
regularity properties for solutions on the entire space $\Bbb{R}^N$, where
the potential $V=V_1+V_2$ belongs to the Kato-Rellich class (i.e., $V_1\in
L^p\left( R^N\right) $, $V_2\in L^\infty \left( R^N\right) $ with $p$
depending on $N$). The proof of Theorem\ref{thm2.2} then systematically
addresses two fundamental cases: (i) when $V_1=0$ almost everywhere, and
(ii) when $\left\| V_1\right\| _{L^p}>0$. In each case, starting from the $%
H^1$ regularity of a weak solution, a distinct iterative scheme is
constructed. The iteration replies on Sobolev embeddings for case (i),
whereas for case (ii) it is driven by a nonlinear recurrence relation
derived from the $L^p$ integrability of $V_1$ via H\"older's inequality.
This scheme leverages the Sobolev embeddings $W^{2,q}\left( R^N\right)
\hookrightarrow L^{q^{*}}\left( R^N\right) $ and precise $L^p$-type
estimates derived from the equation to successively boost the solution's
integrability. A crucial feature is the explicit control over the number of
iteration steps: it is finite, determined intrinsically by the dimension $N$
and the exponent $p$(for case (ii)), and independent of the parameter $%
\lambda $. This process not only yields the specific regularity conclusions (%
$C_B^1$,$C_B^0$ or membership in certain $L^q$ spaces as stated) but, most
importantly, provides a uniform bound in $L^\infty \left( R^N\right) $ for
any bounded family of solutions in $H^1$. This uniform $L^\infty \left(
R^N\right) $ estimate is indispensable for the subsequent compactness and
convergence arguments developed in Section 5.

In what follows, for $a_1$,$a_2$,$\cdots $,$a_m\in \Bbb{R}$, $m\geq 2$, we
stipulate the definition $\Lambda \left( a_1\text{,}\cdots \text{,}%
a_m\right) :=\max \left\{ a_1,a_2,\cdots ,a_m\right\} $, $\Gamma \left( a_1%
\text{,}\cdots \text{,}a_m\right) :=\min \left\{ a_1,a_2,\cdots ,a_m\right\}
$. We first take a look back at the celebrated Agmon-Douglis-Nirenberg
theorem with Dirichlet boundary condition(see \cite{ADN1},\cite{ADN2}):

\begin{proposition}
\label{prop2.1} Let $\Omega $ be a bounded domain in $\Bbb{R}^N$ with smooth
boundary $\partial \Omega \in C^{2+\alpha }$, and $L$ be a uniformly
elliptic operator on $\Omega $, s.t.
\begin{eqnarray*}
Lu=-\sum\limits_{i,j=1}^na_{ij}\left( x\right) \frac{\partial ^2u\left(
x\right) }{\partial x_i\partial x_j}+\sum\limits_{i=1}^nb_i\left( x\right)
\frac{\partial u\left( x\right) }{\partial x_i},\text{\textrm{\ }}x\in
\Omega , \\
u=0,\text{\textrm{\ }}x\in \partial \Omega ,
\end{eqnarray*}
$a_{ij}\left( x\right) $, $b_i\left( x\right) \in C\left( \overline{\Omega }%
\right) $. Consider the following boundary value problem:
\begin{equation}
\left\{
\begin{array}{cc}
Lu+c\left( x\right) u=f\left( x\right) , & x\in \Omega , \\
u=0, & x\in \partial \Omega ,
\end{array}
\right.  \label{eq5}
\end{equation}
and assume $a_{ij}\left( x\right) $, $b_i\left( x\right) $, $c\left(
x\right) \in C^\alpha \left( \overline{\Omega }\right) $, $c\left( x\right)
\geq 0$. Suppose $f\left( x\right) $ $\in C^\alpha \left( \overline{\Omega }%
\right) $, then $\left( \ref{eq5}\right) $ admits a unique solution $u\left(
x\right) \in C^{2+\alpha }\left( \overline{\Omega }\right) $ with Schauder
estimate
\begin{equation}
\left\| u\right\| _{C^{2+\alpha }\left( \overline{\Omega }\right) }\leq
C_1\left\| f\right\| _{C^\alpha \left( \overline{\Omega }\right) },
\label{eq6}
\end{equation}
and if $f\left( x\right) $ $\in L^p\left( \Omega \right) $, $p>1$, then $%
\left( \ref{eq5}\right) $ possesses a unique solution $u\left( x\right) \in
W^{2,p}\left( \Omega \right) $ with $L_p$ estimate
\begin{equation}
\left\| u\right\| _{W^{2,p}\left( \Omega \right) }\leq C_2\left\| f\right\|
_{L^p\left( \Omega \right) },  \label{eq7}
\end{equation}
where $C_1$,$C_2$ are independent of $u$ and $f$.
\end{proposition}

\textbf{A natural question which arises is: Can we establish the
corresponding Schauder estimate or }$L_p$ \textbf{estimate for the case }$%
c\left( x\right) \geq 0$\textbf{,} $c\left( x\right) \neq 0$\textbf{, }$%
c\left( x\right) \in L^\infty \left( \Omega \right) $ \textbf{or }$c\left(
x\right) \in L^p\left( \Omega \right) $\textbf{,} $p>1$\textbf{?}

Aiming at answering this above-mentioned question, we first investigate the
global regularity of solutions of Schr\"odinger equations:

\begin{theorem}
\label{thm2.2} Given $\lambda \in \Bbb{R}$ and let $V=V_1+V_2$ be a K-R
potential, i.e., $V_1\in L^p$, $V_2$ $\in L^\infty \left( \Bbb{R}^N\right) $%
, with $p=2$ if $N\leq 3$, and $p>\frac N2$ if $N\geq 4$. Assume $g\in
C\left( \Bbb{R},\Bbb{R}\right) $ and there exists a $C>0$, $\left| g\left(
t\right) \right| \leq C\left| t\right| $ for $t\in \Bbb{R}$. Suppose $u$ is
a weak nontrivial solution of $\left( \ref{eq1}\right) $.

$\left( i\right) $ Let $V_1=0$ a.e. on $\Bbb{R}^N$. Then $u\in C_B^1\left(
\Bbb{R}^N\right) $;

$\left( ii\right) $ Let $\left\| V_1\right\| _{L^p}>0$. If $p>N\geq 4$ or $%
p=2$ and $N=1$, $u\in C_B^1\left( \Bbb{R}^N\right) $, and if $\frac N2<p\leq
N$ with $N\geq 4$ or $p=2$ and $N=3$, $u\in C_B^0\left( \Bbb{R}^N\right) $,
and especially $p=N=2$, $u\in L^q\left( \Bbb{R}^N\right) $ for $\forall q\in
\left[ 2,\infty \right) $, where
\[
C_B^j\left( \Bbb{R}^N\right) :=\left\{ u\in C^j\left( \Bbb{R}^N\right) \text{%
, }\nabla ^\alpha u\text{ is bounded on }\Bbb{R}^N\text{, }\left| \alpha
\right| \leq j\right\}
\]
for $j\in \Bbb{N\cup }\left\{ 0\right\} $(Concerning the definition of $%
C_B^j\left( \Omega \right) $ for a more general domain $\Omega $, the
readers may consult \cite{AF} for details).
\end{theorem}

\textit{Proof. }Denote by $\rho \left( A\right) $ the resolvent set of $A$.
Take $\widetilde{\lambda }\in \left( -\infty ,-m\right) \cap \rho \left(
A\right) $, s.t. $\widetilde{\lambda }<-m$. Hence $D\left( \widetilde{%
\lambda }I-A\right) =H^2\left( \Bbb{R}^N\right) $ and
\[
R\left( \widetilde{\lambda }I-A\right) =\left( \widetilde{\lambda }%
I-A\right) H^2\left( \Bbb{R}^N\right) =L^2\left( \Bbb{R}^N\right) ,
\]
i.e., for $\forall f\in L^2\left( \Bbb{R}^N\right) $, $\exists \xi \in
H^2\left( \Bbb{R}^N\right) $, s.t.
\[
\left( \widetilde{\lambda }I-A\right) \xi =f.
\]
Due to $g\left( u\right) \in L^2\left( \Bbb{R}^N\right) $, there exists a
unique $\xi ^{*}\in H^2\left( \Bbb{R}^N\right) $, s.t.
\begin{equation}
\left( \widetilde{\lambda }I-A\right) \xi ^{*}=-\left( \lambda -\widetilde{%
\lambda }\right) u-g\left( u\right) .  \label{eq8}
\end{equation}
Observe that $\left( \ref{eq1}\right) $ is equivalent to
\begin{equation}
\widetilde{\lambda }\int_{\Bbb{R}^N}u\varphi -\int_{\Bbb{R}^N}\nabla u\nabla
\varphi +Vu\varphi =-\int_{\Bbb{R}^N}\left[ \left( \lambda -\widetilde{%
\lambda }\right) u+g\left( u\right) \right] \varphi ,  \label{eq9}
\end{equation}
for\textrm{\ }$\forall \varphi \in H^1\left( \Bbb{R}^N\right) $. Combining $%
\left( \ref{eq8}\right) $ and $\left( \ref{eq9}\right) $,
\[
\widetilde{\lambda }\int_{\Bbb{R}^N}\left( \xi ^{*}-u\right) \varphi -\nabla
\left( \xi ^{*}-u\right) \nabla \varphi -V\left( \xi ^{*}-u\right) \varphi
=0,
\]
yielding
\begin{equation}
\left\| \xi ^{*}-u\right\| _X^2=\left( \widetilde{\lambda }+m\right) \left\|
\xi ^{*}-u\right\| _{L^2}^2\leq 0  \label{eq10}
\end{equation}
and this alludes to $u\in H^2\left( \Bbb{R}^N\right) $ via $\left( \ref{eq7}%
\right) $.

Two cases are taken into account:

$\left( i\right) $ $V_1=0$ a.e. on $\Bbb{R}^N$. Divide the proof into two
parts. The foremost task is to verify $u\in L^\infty \left( \Bbb{R}^N\right)
$. Let $p_0=2$. Then $u\in H^2\left( \Bbb{R}^N\right) =W^{2,2}\left( \Bbb{R}%
^N\right) \hookrightarrow L^\infty \left( \Bbb{R}^N\right) $ if $N\leq 3$,
and $u\in W^{2,2}\left( \Bbb{R}^N\right) \hookrightarrow L^q\left( \Bbb{R}%
^N\right) $ for $q\in \left[ 2,+\infty \right) $ if $N=4$, and $u\in
W^{2,2}\left( \Bbb{R}^N\right) \hookrightarrow L^{\frac{2N}{N-4}}\left( \Bbb{%
R}^N\right) $ if $N\geq 5$. Plainly $u\in W^{2,q}\left( \Bbb{R}^N\right)
\hookrightarrow L^\infty \left( \Bbb{R}^N\right) $ for $q>2$ and $N=4$, so
we focus on the cases $N\geq 5$.

Let $p_1=\frac{p_0N}{N-2p_0}$. Clearly $u\in W^{2,p_1}\left( \Bbb{R}%
^N\right) \hookrightarrow L^\infty \left( \Bbb{R}^N\right) $ as $5\leq N\leq
7$. If $N=8$, then $u\in W^{2,p_1}\left( \Bbb{R}^N\right) \hookrightarrow
L^q\left( \Bbb{R}^N\right) $ for $q\in \left[ p_1,+\infty \right) $ via $%
p_1=\frac N2$. Take $q\in \left( p_1,+\infty \right) $ and we get
\begin{equation}
\int_{\Bbb{R}^N}\left| \Delta u\right| ^qdx\leq \left( \left\| V_2\right\|
_{L^\infty }+\left| \lambda \right| +C\right) ^q\int_{\Bbb{R}^N}\left|
u\right| ^qdx.  \label{eq11}
\end{equation}
Since for an integer $r\geq 2$ and $p\in \left[ 1,+\infty \right) $, there
exists $K=K\left( r,p,N\right) >0$, for $\forall \varepsilon >0$ and $j\in
\Bbb{N}$, $0\leq j\leq r-1$, $\forall u\in W^{r,p}\left( \Bbb{R}^N\right) $,
\begin{equation}
\left\| u\right\| _{W^{j,p}}\leq K\varepsilon \left\| u\right\|
_{W^{r,p}}+K\varepsilon ^{-\frac j{r-j}}\left\| u\right\| _{L^p}
\label{eq12}
\end{equation}
(see \cite{AF},\cite{B}). Therefore
\begin{eqnarray}
\left\| \nabla u\right\| _{L^q} &\leq &K\left( 2,q,N\right) \varepsilon
\left\| u\right\| _{W^{2,q}}+K\left( 2,q,N\right) \varepsilon ^{-1}\left\|
u\right\| _{L^q}  \nonumber \\
\ &\leq &K\left( 2,q,N\right) \varepsilon \left( \left\| \Delta u\right\|
_{L^q}+\left\| \nabla u\right\| _{L^q}+\left\| u\right\| _{L^q}\right)
+K\left( 2,q,N\right) \varepsilon ^{-1}\left\| u\right\| _{L^q}.
\label{eq13}
\end{eqnarray}
Take $\varepsilon \leq \frac 1{2K\left( 2,q,N\right) }$ and we can easily
infer that
\begin{equation}
\left\| \nabla u\right\| _{L^q}\leq \left\| \Delta u\right\| _{L^q}+\left[
1+2K\left( 2,q,N\right) \varepsilon ^{-1}\right] \cdot \left\| u\right\|
_{L^q}.  \label{eq14}
\end{equation}
Inserting $\left( \ref{eq14}\right) $ into $\left( \ref{eq11}\right) $,
\begin{equation}
\left\| \nabla u\right\| _{L^q}\leq \left[ \left( \left\| V_2\right\|
_{L^\infty }+\left| \lambda \right| +C\right) +2K\left( 2,q,N\right)
\varepsilon ^{-1}+1\right] \cdot \left\| u\right\| _{L^q}.  \label{eq15}
\end{equation}

Notice that for for an integer $r\geq 1$ and $s\in \left[ 1,+\infty \right) $%
, if $s<\frac Nr$, there exists $C\left( s,N\right) >0$,
\[
C\left( r,s,N\right) \left\| u\right\| _{W^{r,s}}\geq \left\| u\right\|
_{L^{s^{*}\left( r\right) }}
\]
for $s^{*}\left( r\right) :=\frac{sN}{N-rs}$, and if $s=\frac Nr$, there
exists $\widetilde{C}\left( r,s,t,N\right) >0$,
\[
\widetilde{C}\left( r,s,t,N\right) \left\| u\right\| _{W^{r,s}}\geq \left\|
u\right\| _{L^t}
\]
for $t\in \left[ s,+\infty \right) $, and if $s>\frac Nr$, there exists $%
C_0\left( s,N\right) >0$,
\[
C_0\left( r,s,N\right) \left\| u\right\| _{W^{r,s}}\geq \left\| u\right\|
_{L^\infty }.
\]
Consequently
\begin{eqnarray*}
\left\| u\right\| _{L^\infty } &\leq &C_0\left( 2,q,N\right) \left\|
u\right\| _{W^{2,q}} \\
\ &\leq &C_0\left( 2,q,N\right) \cdot 2\left( \left\| V_2\right\| _{L^\infty
}+\left| \lambda \right| +C+K\left( 2,q,N\right) \varepsilon ^{-1}+1\right)
\cdot \left\| u\right\| _{L^q}
\end{eqnarray*}
via $\left( \ref{eq11}\right) $ and $\left( \ref{eq15}\right) $. Denote
\[
\gamma \left( q,\lambda ,N,\varepsilon \right) :=2\left( \left\| V_2\right\|
_{L^\infty }+\left| \lambda \right| +C+K\left( 2,q,N\right) \varepsilon
^{-1}+1\right) .
\]
Hence
\begin{eqnarray}
\left\| u\right\| _{L^\infty } &\leq &C_0\left( 2,q,N\right) \widetilde{C}%
\left( 2,p_1,q,N\right) C\left( 2,p_0,N\right) \gamma \left( q,\lambda
,N,\varepsilon \right) \gamma \left( p_1,\lambda ,N,\varepsilon \right)
\gamma \left( p_0,\lambda ,N,\varepsilon \right) \left\| u\right\| _{H^1}.
\nonumber \\
&&  \label{eq16}
\end{eqnarray}

Let $N\geq 9$, $p_2=\frac{p_1N}{N-2p_1}$,$\cdots $, $p_j=\frac{p_{j-1}N}{%
N-2p_{j-1}}$, and assume $j\geq 1$, $p_i<\frac N2$, $p_{j+1}\geq \frac N2$, $%
i=0$,$\cdots $,$j$. It is pretty straightforward for us to show that $p_i=%
\frac{2N}{N-2p_0i}$ by induction as $i=1$,$\cdots $,$j+1$. Observe that
\begin{eqnarray*}
p_j &<&\frac N2\Leftrightarrow j<\frac N4-1; \\
p_{j+1} &\geq &\frac N2\Leftrightarrow j\geq \frac N4-2.
\end{eqnarray*}
Hence $j=\frac N4-2$ if $N=2k$ and $k$ even, $j=\frac N4-\frac 32$ if $N=2k$
and $k$ odd, $j=\frac{N-5}4$ if $N=2k+1$ and $k$ even, $j=\frac{N-7}4$ if $%
N=2k+1$ and $k$ odd, which explicitly manifests that $p_j=\frac N4$ and $%
p_{j+1}=\frac N2$ if $N=2k$ and $k$ even, $p_j=\frac N3$ and $p_{j+1}=N$ if $%
N=2k$ and $k$ odd, $p_j=\frac{2N}5$ and $p_{j+1}=2N$ if $N=2k+1$ and $k$
even, $p_j=\frac{2N}7$ and $p_{j+1}=\frac{2N}3$ if $N=2k+1$ and $k$ odd.

Take $\varepsilon <\Gamma \left( \frac 1{2K\left( 2,p_0,N\right) },\frac
1{2K\left( 2,p_1,N\right) },\cdots ,\frac 1{2K\left( 2,p_{j+1},N\right)
}\right) $. By an alternative bootstrap iteration argument, we obtain
\begin{eqnarray}
\left\| u\right\| _{L^\infty } &\leq &C_0\left( 2,p_{j+1},N\right) \left\|
u\right\| _{W^{2,p_{j+1}}}  \nonumber \\
\ &\leq &C_0\left( 2,p_{j+1},N\right) \gamma \left( p_{j+1},\lambda
,N,\varepsilon \right) \left\| u\right\| _{L^{p_{j+1}}}  \nonumber \\
\ &\leq &C_0\left( 2,p_{j+1},N\right) \gamma \left( p_{j+1},\lambda
,N,\varepsilon \right) C\left( 2,p_j,N\right) \left\| u\right\| _{W^{2,p_j}}
\nonumber \\
\ &\leq &C_0\left( 2,p_{j+1},N\right) C\left( 2,p_j,N\right) \gamma \left(
p_{j+1},\lambda ,N,\varepsilon \right) \gamma \left( p_j,\lambda
,N,\varepsilon \right) \left\| u\right\| _{L^{p_j}}  \nonumber \\
\ &\leq &C_0\left( 2,p_{j+1},N\right) \cdot \prod\limits_{i=0}^jC\left(
2,p_i,N\right) \cdot \prod\limits_{i=0}^{j+1}\gamma \left( p_i,\lambda
,N,\varepsilon \right) \cdot \left\| u\right\| _{H^1}  \label{eq17}
\end{eqnarray}
if $N$ is odd or $N=2k$ with $k$ odd.

If $N=2k$ with $k$ even, $W^{2,p_{j+1}}\left( \Bbb{R}^N\right)
\hookrightarrow L^q\left( \Bbb{R}^N\right) $ for $q\in \left[
p_{j+1},+\infty \right) $ due to $p_{j+1}=\frac N2$. Choose $q\in \left(
p_{j+1},+\infty \right) $ and pick $\varepsilon <\Gamma \left( \frac
1{2K\left( 2,p_0,N\right) },\frac 1{2K\left( 2,p_1,N\right) },\cdots ,\frac
1{2K\left( 2,p_{j+1},N\right) },\frac 1{2K\left( 2,q,N\right) }\right) $.
Thus
\begin{eqnarray}
\left\| u\right\| _{L^\infty } &\leq &C_0\left( 2,q,N\right) \left\|
u\right\| _{W^{2,q}}  \nonumber \\
\ &\leq &C_0\left( 2,q,N\right) \gamma \left( q,\lambda ,N,\varepsilon
\right) \left\| u\right\| _{L^q}  \nonumber \\
\ &\leq &C_0\left( 2,q,N\right) \gamma \left( q,\lambda ,N,\varepsilon
\right) \widetilde{C}\left( 2,p_{j+1},q,N\right) \left\| u\right\|
_{W^{2,p_{j+1}}}  \nonumber \\
\ &\leq &C_0\left( 2,q,N\right) \widetilde{C}\left( 2,p_{j+1},q,N\right)
\cdot \prod\limits_{i=0}^jC\left( 2,p_i,N\right) \cdot \gamma \left(
q,\lambda ,N,\varepsilon \right) \cdot \prod\limits_{i=0}^{j+1}\gamma \left(
p_i,\lambda ,N,\varepsilon \right) \cdot \left\| u\right\| _{H^1}.  \nonumber
\\
&&  \label{eq18}
\end{eqnarray}
We are left to confirm the validity of $u\in C_B^1\left( \Bbb{R}^N\right) $
and initially deal with the case $N\geq 4$. Let $\widetilde{q}>N$ and
observe that $N\geq 2^{*}=\frac{2N}{N-2}\Leftrightarrow N\geq 4$. Hence
\begin{eqnarray}
\left\| \Delta u\right\| _{L^{\widetilde{q}}} &\leq &\left( \left\|
V_2\right\| _{L^\infty }+\left| \lambda \right| +C\right) \left\| u\right\|
_{L^{\widetilde{q}}}  \nonumber \\
\ &\leq &\left( \left\| V_2\right\| _{L^\infty }+\left| \lambda \right|
+C\right) \left\| u\right\| _{L^\infty }^{1-\frac{2^{*}}{\widetilde{q}}%
}\cdot \left\| u\right\| _{L^{2^{*}}}^{\frac{2^{*}}{\widetilde{q}}}
\nonumber \\
\ &\leq &\left( \left\| V_2\right\| _{L^\infty }+\left| \lambda \right|
+C\right) \left\| u\right\| _{L^\infty }^{1-\frac{2^{*}}{\widetilde{q}}%
}L\left( 2^{*},N\right) ^{\frac{2^{*}}{\widetilde{q}}}\left\| u\right\|
_{H^1}^{\frac{2^{*}}{\widetilde{q}}},  \label{eq19}
\end{eqnarray}
where $L\left( 2^{*},N\right) $ is the embedding constant arising from $%
\left( \ref{eq28}\right) $.

Take $\varepsilon \leq \frac 1{2K\left( 2,\widetilde{q},N\right) }$ and
therefore
\[
\left\| \nabla u\right\| _{L^{\widetilde{q}}}\leq \left\| \Delta u\right\|
_{L^{\widetilde{q}}}+\left[ 1+2K\left( 2,\widetilde{q},N\right) \varepsilon
^{-1}\right] \cdot \left\| u\right\| _{L^{\widetilde{q}}},
\]
deriving
\begin{eqnarray*}
\left\| u\right\| _{W^{2,\widetilde{q}}} &\leq &\gamma \left( \widetilde{q}%
,\lambda ,N,\varepsilon \right) \left\| u\right\| _{L^{\widetilde{q}}} \\
&\leq &\gamma \left( \widetilde{q},\lambda ,N,\varepsilon \right) \left\|
u\right\| _{L^\infty }^{1-\frac{2^{*}}{\widetilde{q}}}L\left( 2^{*},N\right)
^{\frac{2^{*}}{\widetilde{q}}}\left\| u\right\| _{H^1}^{\frac{2^{*}}{%
\widetilde{q}}}
\end{eqnarray*}
with
\[
\gamma \left( \widetilde{q},\lambda ,N,\varepsilon \right) =2\left( \left\|
V_2\right\| _{L^\infty }+\left| \lambda \right| +C+K\left( 2,\widetilde{q}%
,N\right) \varepsilon ^{-1}+1\right) .
\]
Being aware that for $r\geq 1$, if $q>\frac Nr$, $W^{r,q}\left( \Bbb{R}%
^N\right) \hookrightarrow L^\infty \left( \Bbb{R}^N\right) $, and moreover
if $r-\frac Nq$ is not an integer and set $k=\left[ r-\frac Nq\right] $,
then for $\forall u\in W^{r,q}\left( \Bbb{R}^N\right) $,
\begin{equation}
\left\| \nabla ^\alpha u\right\| _{L^\infty \left( \Bbb{R}^N\right) }\leq
\sigma \left( q,N\right) \left\| u\right\| _{W^{r,q}},\mathrm{\ }\forall
\alpha ,\mathrm{\ }\left| \alpha \right| \leq k,  \label{eq20}
\end{equation}
and furthermore
\begin{equation}
\left| \nabla ^\alpha u\left( x\right) -\nabla ^\alpha u\left( y\right)
\right| \leq \sigma \left( q,N\right) \left\| u\right\| _{W^{r,q}}\left|
x-y\right| \text{,\textrm{\ }}\forall x,y\in \Bbb{R}^N,\text{ }\forall
\alpha ,\text{ }\left| \alpha \right| <k,\text{ }  \label{eq21}
\end{equation}
and
\begin{equation}
\left| \nabla ^\alpha u\left( x\right) -\nabla ^\alpha u\left( y\right)
\right| \leq \sigma \left( q,N\right) \left\| u\right\| _{W^{r,q}}\left|
x-y\right| ^\theta ,\text{\textrm{\ }a.e. }x,y\in \Bbb{R}^N,\text{\textrm{\ }%
}\forall \alpha ,\text{ }\left| \alpha \right| =k,  \label{eq22}
\end{equation}
$\theta =r-\frac Nq-k$, $0<\theta <1$(see \cite{B}), where $\left[ r-\frac
Nq\right] $ stands for the integer part of $r-\frac Nq$. Take $r=2$, $q=%
\widetilde{q}$, and consequently we assert that the proof of $\left(
i\right) $ follows for $N\geq 4$.

As for the case $N\leq 3$, we need merely to deal with the case $\widetilde{q%
}\in \left( \Lambda \left( 2,N\right) ,2^{*}\right) $. A simple observation
exhibits
\begin{eqnarray}
\left\| \Delta u\right\| _{L^{\widetilde{q}}} &\leq &\left( \left\|
V_2\right\| _{L^\infty }+\left| \lambda \right| +C\right) \left\| u\right\|
_{L^{\widetilde{q}}}  \nonumber \\
&\leq &2\left( \left\| V_2\right\| _{L^\infty }+\left| \lambda \right|
+C\right) \cdot \left( 1+L\left( 2,\widetilde{q},N\right) \right) \left\|
u\right\| _{H^1},  \label{eq23}
\end{eqnarray}
where $L\left( 2,\widetilde{q},N\right) $ is the embedding constant arising
from $\left( \ref{eq29}\right) $. A standard argument arrives at the
assertion and we have thus far concluded the proof of $\left( i\right) $.

$\left( ii\right) $\textbf{\ }$\left\| V_1\right\| _{L^p}>0$. For the sake
of convenience, we divide the proof into three cases:

$\left\langle a\right\rangle $ $p>N\geq 4$. Given $q\in \left( 0,N\right) $,
it is easy to check
\begin{equation}
\frac{qpN}{\left( N-q\right) p+qN}>q\Leftrightarrow p>N,  \label{eq24}
\end{equation}
which indeed is frequently used in the iteration argument.

Let $q_0=2$, $p=q\alpha $, $q\in \left[ q_0,\frac{q_0N}{N-q_0}\right) $.
Clearly $p>q$ is due to
\begin{equation}
N\geq \frac{q_0N}{N-q_0}\Leftrightarrow N\geq 4.  \label{eq25}
\end{equation}
Hence
\begin{equation}
q\beta \leq \frac{q_0N}{N-q_0}\Leftrightarrow q\leq \frac{q_0pN}{\left(
N-q_0\right) p+q_0N},  \label{eq26}
\end{equation}
$\beta =\frac p{p-q}$, $\frac 1\alpha +\frac 1\beta =1$. Set $q_1=\frac{q_0pN%
}{\left( N-q_0\right) p+q_0N}$ and thus $q_1>q_0$ via $\left( \ref{eq24}%
\right) $. Therefore
\begin{eqnarray}
\int_{\Bbb{R}^N}\left| \Delta u\right| ^{q_1}dx &\leq &2^{q_1-1}\left[ \int_{%
\Bbb{R}^N}\left| V_1u\right| ^{q_1}dx+\left( \left\| V_2\right\| _{L^\infty
}+\left| \lambda \right| +C\right) ^{q_1}\int_{\Bbb{R}^N}\left| u\right|
^{q_1}dx\right]  \nonumber \\
&\leq &2^{q_1-1}\left( \int_{\Bbb{R}^N}\left| V_1\right| ^pdx\right) ^{\frac
1{\alpha _1}}\cdot \left( \int_{\Bbb{R}^N}\left| u\right| ^{q_1\beta
_1}dx\right) ^{\frac 1{\beta _1}}  \nonumber \\
&&+2^{q_1-1}\left( \left\| V_2\right\| _{L^\infty }+\left| \lambda \right|
+C\right) ^{q_1}\int_{\Bbb{R}^N}\left| u\right| ^{q_1}dx,  \label{eq27}
\end{eqnarray}
$\alpha _1=\frac p{q_1}=1+\frac{\left( N-q_0\right) p}{q_0N}$, $\beta _1=1+%
\frac{q_0N}{\left( N-q_0\right) p}$. Notice that for $s\in \left[ 1,N\right)
$, there exists $L\left( s^{*},N\right) >0$,
\begin{equation}
L\left( s^{*},N\right) \left\| u\right\| _{W^{1,s}}\geq \left\| u\right\|
_{L^{s^{*}}}  \label{eq28}
\end{equation}
for $s^{*}=\frac{sN}{N-s}$, and for $t\in \left[ s,s^{*}\right] $, there
exists $L\left( s,t,N\right) >0$,
\begin{equation}
L\left( s,t,N\right) \left\| u\right\| _{W^{1,s}}\geq \left\| u\right\|
_{L^t}.  \label{eq29}
\end{equation}
Thereby, $\left( \ref{eq27}\right) $ yields
\begin{eqnarray}
\left\| \Delta u\right\| _{L^{q_1}} &\leq &2^{1-\frac 1{q_1}}\left[ \left\|
V_1\right\| _{L^p}L\left( q_0^{*},N\right) +\left( \left\| V_2\right\|
_{L^\infty }+\left| \lambda \right| +C\right) L\left( 2,q_1,N\right) \right]
\cdot \left\| u\right\| _{H^1}.  \nonumber \\
&&  \label{eq30}
\end{eqnarray}
Define
\[
\rho \left( q_0,q_1,\lambda ,N\right) :=2^{1-\frac 1{q_1}}\left[ \left\|
V_1\right\| _{L^p}L\left( q_0^{*},N\right) +\left( \left\| V_2\right\|
_{L^\infty }+\left| \lambda \right| +C\right) L\left( 2,q_1,N\right) \right]
\]
and thus
\begin{eqnarray}
\left\| \nabla u\right\| _{L^{q_1}} &\leq &\left\| \Delta u\right\|
_{L^{q_1}}+\left[ 1+2K\left( 2,q_1,N\right) \varepsilon ^{-1}\right] \cdot
\left\| u\right\| _{L^{q_1}}  \nonumber \\
&\leq &\left( \rho \left( q_0,q_1,\lambda ,N\right) +\left[ 1+2K\left(
2,q_1,N\right) \varepsilon ^{-1}\right] L\left( 2,q_1,N\right) \right)
\left\| u\right\| _{H^1}.  \label{eq31}
\end{eqnarray}
by taking $\varepsilon \leq \frac 1{2K\left( 2,q_1,N\right) }$. Set
\[
\rho \left( q_0,q_1,\lambda ,N,\varepsilon \right) :=2\rho \left(
q_0,q_1,\lambda ,N\right) +2\left( K\left( 2,q_1,N\right) \varepsilon
^{-1}+1\right) L\left( 2,q_1,N\right) .
\]
Consequently $\left( \ref{eq30}\right) $ along with $\left( \ref{eq31}%
\right) $ derives
\begin{equation}
\left\| u\right\| _{W^{2,q_1}}\leq \rho \left( q_0,q_1,\lambda
,N,\varepsilon \right) \cdot \left\| u\right\| _{H^1}.  \label{eq32}
\end{equation}

Set $p=q\alpha $, $q\in \left[ q_1,q_1^{*}\right) $, $q_1^{*}=\frac{q_1N}{%
N-q_1}$. Thus
\begin{equation}
q\beta \leq q_1^{*}\Leftrightarrow q\leq \frac{q_1pN}{\left( N-q_1\right)
p+q_1N},  \label{eq33}
\end{equation}
$\beta =\frac p{p-q}$, $\frac 1\alpha +\frac 1\beta =1$. Let $q_2=\frac{q_1pN%
}{\left( N-q_1\right) p+q_1N}$. Due to
\begin{equation}
\left( N-q_1\right) p+q_1N>0\Leftrightarrow p\left( N-2q_0\right) +2q_0N>0
\label{eq34}
\end{equation}
which is certainly true for $N\geq 4$, we get $q_2>q_1$ by $\left( \ref{eq24}%
\right) $. Hence
\begin{eqnarray}
\int_{\Bbb{R}^N}\left| \Delta u\right| ^{q_2}dx &\leq &2^{q_2-1}\left[ \int_{%
\Bbb{R}^N}\left| V_1u\right| ^{q_2}dx+\left( \left\| V_2\right\| _{L^\infty
}+\left| \lambda \right| +C\right) ^{q_2}\int_{\Bbb{R}^N}\left| u\right|
^{q_2}dx\right]  \nonumber \\
\ &\leq &2^{q_2-1}\left( \int_{\Bbb{R}^N}\left| V_1\right| ^pdx\right)
^{\frac 1{\alpha _2}}\cdot \left( \int_{\Bbb{R}^N}\left| u\right| ^{q_2\beta
_2}dx\right) ^{\frac 1{\beta _2}}  \nonumber \\
&&\ \ +2^{q_2-1}\left( \left\| V_2\right\| _{L^\infty }+\left| \lambda
\right| +C\right) ^{q_2}\int_{\Bbb{R}^N}\left| u\right| ^{q_2}dx,
\label{eq35}
\end{eqnarray}
yielding
\begin{eqnarray}
\left\| \Delta u\right\| _{L^{q_2}} &\leq &2^{1-\frac 1{q_2}}\left\|
V_1\right\| _{L^p}\cdot \left\| u\right\| _{L^{q_1^{*}}}+2^{1-\frac
1{q_2}}\left( \left\| V_2\right\| _{L^\infty }+\left| \lambda \right|
+C\right) \cdot \left\| u\right\| _{L^{q_2}}  \nonumber \\
\ &\leq &2^{1-\frac 1{q_2}}\left[ \left\| V_1\right\| _{L^p}L\left(
q_1^{*},N\right) +\left( \left\| V_2\right\| _{L^\infty }+\left| \lambda
\right| +C\right) L\left( q_1,q_2,N\right) \right] \cdot \left\| u\right\|
_{W^{1,q_1}},  \nonumber \\
&&  \label{eq36}
\end{eqnarray}
$\alpha _2=\frac p{q_2}=1+\frac{\left( N-q_1\right) p}{q_1N}$, $\beta _2=1+%
\frac{q_1N}{\left( N-q_1\right) p}$. Set
\[
\rho \left( q_1,q_2,\lambda ,N\right) =2^{1-\frac 1{q_2}}\left[ \left\|
V_1\right\| _{L^p}L\left( q_1^{*},N\right) +\left( \left\| V_2\right\|
_{L^\infty }+\left| \lambda \right| +C\right) L\left( q_1,q_2,N\right)
\right]
\]
and thus
\begin{eqnarray}
\left\| \nabla u\right\| _{L^{q_2}} &\leq &\left\| \Delta u\right\|
_{L^{q_2}}+\left[ 1+2K\left( 2,q_2,N\right) \varepsilon ^{-1}\right] \cdot
\left\| u\right\| _{L^{q_2}}  \nonumber \\
\ &\leq &\left( \rho \left( q_1,q_2,\lambda ,N\right) +\left[ 1+2K\left(
2,q_1,N\right) \varepsilon ^{-1}\right] L\left( q_1,q_2,N\right) \right)
\left\| u\right\| _{W^{1,q_1}}  \label{eq37}
\end{eqnarray}
by taking $\varepsilon \leq \frac 1{2K\left( 2,q_2,N\right) }$. Set
\[
\rho \left( q_1,q_2,\lambda ,N,\varepsilon \right) :=2\rho \left(
q_1,q_2,\lambda ,N\right) +2\left( K\left( 2,q_1,N\right) \varepsilon
^{-1}+1\right) L\left( q_1,q_2,N\right) .
\]
Thereby, $\left( \ref{eq36}\right) $ coupling with $\left( \ref{eq37}\right)
$ gives rise to
\begin{equation}
\left\| u\right\| _{W^{2,q_2}}\leq \rho \left( q_1,q_2,\lambda
,N,\varepsilon \right) \cdot \left\| u\right\| _{W^{1,q_1}}.  \label{eq38}
\end{equation}

Let $q_3=\frac{q_2pN}{\left( N-q_2\right) p+q_2N}$, $\cdots $, $q_j=\frac{%
q_{j-1}pN}{\left( N-q_{j-1}\right) p+q_{j-1}N}$, $q_{j+1}=\frac{q_jpN}{%
\left( N-q_j\right) p+q_jN}$, and for the sake of argument we suppose $%
q_i<\frac N2$ for $i=0$,$1$,$\cdots $,$j$. Based on this assumption, the
denotation of $q_i$ undoubtedly makes sense for $i=1$,$\cdots $,$j+1$. And
also assume $q_{j+1}\geq \frac N2$. It is conspicuously deduced by induction
that $q_i=\frac{2pN}{\left( N-2i\right) p+2iN}$ for $i=0$,$1$,$\cdots $,$j+1$%
. Notice that $N\geq 2^{*}$ as $N\geq 4$. In terms of the hypothesis,
\begin{eqnarray}
q_j &<&\frac N2\Leftrightarrow j<\frac{\left( N-4\right) p}{2\left(
p-N\right) };  \label{eq39} \\
q_{j+1} &\geq &\frac N2\Leftrightarrow j\geq \frac{\left( N-4\right) p}{%
2\left( p-N\right) }-1.  \label{eq40}
\end{eqnarray}
Objectively, for $N\geq 5$, we can find such a $j_0\in \Bbb{N\cup }\left\{
0\right\} $, s.t., $q_{j_0}<\frac N2$, $q_{j_0+1}\geq \frac N2$. By way of
negation, we have
\begin{eqnarray*}
q_i &=&\frac{2pN}{\left( N-2i\right) p+2iN} \\
&=&\frac{2pN}{Np-2i\left( p-N\right) }<0
\end{eqnarray*}
as $i>\frac{Np}{2\left( P-N\right) }$, subverting the assumption!

Define
\begin{eqnarray*}
\Theta ^{*} &=&\left\{ i\in \Bbb{N\cup }\left\{ 0\right\} :i<\frac{\left(
N-4\right) p}{2\left( p-N\right) }\right\} ; \\
\Theta _{*} &=&\left\{ i\in \Bbb{N\cup }\left\{ 0\right\} :i\geq \frac{%
\left( N-4\right) p}{2\left( p-N\right) }-1\right\} .
\end{eqnarray*}
Clearly $\Theta ^{*}=\varnothing $ and $\inf \Theta _{*}=0$ as $N=4$. Notice
that
\begin{eqnarray}
\frac{\left( N-4\right) p}{2\left( p-N\right) } &>&\frac 12\Leftrightarrow
\left( N-5\right) p>-N;  \label{eq41} \\
\frac{\left( N-4\right) p}{2\left( p-N\right) } &>&1\Leftrightarrow \left(
N-6\right) p>-2N,  \label{eq42}
\end{eqnarray}
so if $N=5$ with $p\geq 2N$, we get $\frac 12<\frac{\left( N-4\right) p}{%
2\left( p-N\right) }\leq 1$ and $\Theta ^{*}\cap \Theta _{*}=\left\{
0\right\} $ accordingly, and if $N=5$ with $\frac 43N\leq p<2N$, we yield $%
\Theta ^{*}\cap \Theta _{*}=\left\{ 1\right\} $ due to $\frac{\left(
N-4\right) p}{2\left( p-N\right) }\leq 2\Leftrightarrow p\geq \frac 43N$.
More generally, one can easily infer that $j_0=l$ if $\frac{\left(
N-4\right) p}{2\left( p-N\right) }\in \left( l,l+1\right) $, and $j_0=l-1$
if $\frac{\left( N-4\right) p}{2\left( p-N\right) }=l$, and here we denote
by $\left[ \frac{\left( N-4\right) p}{2\left( p-N\right) }\right] $ the
integer part of $\frac{\left( N-4\right) p}{2\left( p-N\right) }$ and set $%
l=\left[ \frac{\left( N-4\right) p}{2\left( p-N\right) }\right] $. Therefore
$j_0\geq 1$ if $N\geq 6$, or $N=5$ with $p<2N$, and $j_0=0$ if $N=5$ with $%
p\geq 2N$.

Take $\varepsilon <\Gamma \left( \frac 1{2K\left( 2,q_0,N\right) },\frac
1{2K\left( 2,q_1,N\right) },\cdots ,\frac 1{2K\left( 2,q_{j_0+1},N\right)
}\right) $. If $N\geq 6$, or $N=5$ with $p<2N$, and moreover $%
q_{j_0+1}>\frac N2$, then a standard iteration argument derives
\begin{eqnarray}
\left\| u\right\| _{L^\infty } &\leq &C_0\left( 2,q_{j_0+1},N\right) \left\|
u\right\| _{W^{2,q_{j_0+1}}}  \nonumber \\
\ &\leq &C_0\left( 2,q_{j_0+1},N\right) \rho \left(
q_{j_0},q_{j_0+1},\lambda ,N,\varepsilon \right) \left\| u\right\|
_{W^{1,q_{j_0}}}  \nonumber \\
\ &\leq &C_0\left( 2,q_{j_0+1},N\right) \rho \left(
q_{j_0},q_{j_0+1},\lambda ,N,\varepsilon \right) \left\| u\right\|
_{W^{2,q_{j_0}}}  \nonumber \\
\ &\leq &C_0\left( 2,q_{j_0+1},N\right) \rho \left(
q_{j_0},q_{j_0+1},\lambda ,N,\varepsilon \right) \rho \left(
q_{j_0-1},q_{j_0},\lambda ,N,\varepsilon \right) \cdot \left\| u\right\|
_{W^{1,q_{j_0-1}}}  \nonumber \\
\ &\leq &C_0\left( 2,q_{j_0+1},N\right) \prod\limits_{i=0}^{j_0}\rho \left(
q_i,q_{i+1},\lambda ,N,\varepsilon \right) \left\| u\right\| _{H^1}.
\label{eq43}
\end{eqnarray}
And if $j_0=\frac{\left( N-4\right) p}{2\left( p-N\right) }-1$, it alludes
to $q_{j_0+1}=\frac N2$. Notice that $W^{2,q_{j_0+1}}\hookrightarrow L^q$
for $\forall q\in \left[ q_{j_0+1},+\infty \right) $. By choosing $q\in
\left( q_{j_0+1},p\right) $,
\begin{eqnarray*}
\int_{\Bbb{R}^N}\left| \Delta u\right| ^qdx &\leq &2^{q-1}\left[ \int_{\Bbb{R%
}^N}\left| V_1u\right| ^qdx+\left( \left\| V_2\right\| _{L^\infty }+\left|
\lambda \right| +C\right) ^q\int_{\Bbb{R}^N}\left| u\right| ^qdx\right] \\
\ &\leq &2^{q-1}\left( \int_{\Bbb{R}^N}\left| V_1\right| ^pdx\right) ^{\frac
1\alpha }\cdot \left( \int_{\Bbb{R}^N}\left| u\right| ^{q\beta }dx\right)
^{\frac 1\beta } \\
&&\ \ \ +2^{q-1}\left( \left\| V_2\right\| _{L^\infty }+\left| \lambda
\right| +C\right) ^q\int_{\Bbb{R}^N}\left| u\right| ^qdx,
\end{eqnarray*}
$\alpha =\frac pq$, $\frac 1\alpha +\frac 1\beta =1$. Take
\[
\varepsilon <\Gamma \left( \frac 1{2K\left( 2,q_0,N\right) },\frac
1{2K\left( 2,q_1,N\right) },\cdots ,\frac 1{2K\left( 2,q_{j_0+1},N\right)
},\frac 1{2K\left( 2,q,N\right) }\right) ,
\]
and set
\begin{eqnarray*}
\widetilde{\rho }\left( q_{j_0+1},q,\beta ,\lambda ,N,\varepsilon \right)
&:&=2^{2-\frac 1q}\left\| V_1\right\| _{L^p}\widetilde{C}\left(
2,q_{j_0+1},q\beta ,N\right) \\
&&\ +\left[ 2^{2-\frac 1q}\left( \left\| V_2\right\| _{L^\infty }+\left|
\lambda \right| +C\right) +2K\left( 2,q,N\right) \varepsilon ^{-1}+2\right]
\widetilde{C}\left( 2,q_{j_0+1},q,N\right) ,
\end{eqnarray*}
we get
\begin{eqnarray}
\ \left\| u\right\| _{W^{2,q}} &\leq &2^{2-\frac 1q}\left\| V_1\right\|
_{L^p}\cdot \left\| u\right\| _{L^{q\beta }}  \nonumber \\
&&\ +\left[ 2^{2-\frac 1q}\left( \left\| V_2\right\| _{L^\infty }+\left|
\lambda \right| +C\right) +2K\left( 2,q,N\right) \varepsilon ^{-1}+2\right]
\cdot \left\| u\right\| _{L^q}  \nonumber \\
&\leq &\widetilde{\rho }\left( q_{j_0+1},q,\beta ,\lambda ,N,\varepsilon
\right) \left\| u\right\| _{W^{2,q_{j_0+1}}},  \label{eq44}
\end{eqnarray}
and hence
\begin{eqnarray}
\left\| u\right\| _{L^\infty } &\leq &C_0\left( 2,q,N\right) \left\|
u\right\| _{W^{2,q}}  \nonumber \\
\ &\leq &C_0\left( 2,q,N\right) \widetilde{\rho }\left( q_{j_0+1},q,\beta
,\lambda ,N,\varepsilon \right) \left\| u\right\| _{W^{2,q_{j_0+1}}}
\nonumber \\
\ &\leq &C_0\left( 2,q,N\right) \widetilde{\rho }\left( q_{j_0+1},q,\beta
,\lambda ,N,\varepsilon \right) \prod\limits_{i=0}^{j_0}\rho \left(
q_i,q_{i+1},\lambda ,N,\varepsilon \right) \left\| u\right\| _{H^1}
\label{eq45}
\end{eqnarray}
and subsequently produce
\begin{eqnarray}
\left\| u\right\| _{W^{2,p}} &\leq &2^{2-\frac 1p}\left\| u\right\|
_{L^\infty }\cdot \left\| V_1\right\| _{L^p}+\beta \left( p,\lambda
,N,\varepsilon \right) \cdot \left\| u\right\| _{L^\infty }^{1-\frac{2^{*}}%
p}\cdot \left\| u\right\| _{L^{2^{*}}}^{\frac{2^{*}}p}  \nonumber \\
\ &\leq &2^{2-\frac 1p}\left\| u\right\| _{L^\infty }\cdot \left\|
V_1\right\| _{L^p}+\beta \left( p,\lambda ,N,\varepsilon \right) \cdot
L\left( 2^{*},N\right) ^{\frac{2^{*}}p}\cdot \left\| u\right\| _{L^\infty
}^{1-\frac{2^{*}}p}\cdot \left\| u\right\| _{H^1}^{\frac{2^{*}}p}  \nonumber
\\
&&\ \ \   \label{eq46}
\end{eqnarray}
also by taking $\varepsilon <\frac 1{2K\left( 2,p,N\right) }$, where
\[
\beta \left( p,\lambda ,N,\varepsilon \right) :=2^{2-\frac 1p}\left( \left\|
V_2\right\| _{L^\infty }+\left| \lambda \right| +C\right) +2K\left(
2,p,N\right) \varepsilon ^{-1}+2.
\]
Therefore $u\in C_B^1\left( \Bbb{R}^N\right) $. And the cases $N=4$, and $%
N=5 $ with $p\geq 2N$ are trivial.

$\left\langle b\right\rangle $ $\frac N2<p\leq N$ with $N\geq 4$. It is
worth noticing that in this case $\left( \ref{eq24}\right) $ is plainly
violated for $q\geq 0$, and so the preceding iterative method presented by
case $\left\langle a\right\rangle $ is implicitly invalidated. For this
reason, we introduce a slightly different measure, which does the trick.
Given $q\in \left( 0,\frac N2\right) $, the following assertion
\begin{equation}
q\leq \frac N2\Rightarrow \left( N-2q\right) p+qN>0  \label{eq47}
\end{equation}
is self-evident and consequently
\begin{equation}
\frac{qpN}{\left( N-2q\right) p+qN}>q\Leftrightarrow p>\frac N2.
\label{eq48}
\end{equation}

Let $N\geq 5$, $q_0=2$, $q_1=\frac{q_0pN}{\left( N-2q_0\right) p+q_0N}$. Due
to $q_0<\frac N2$, we obtain $q_1>q_0$ via $\left( \ref{eq47}\right) $ and $%
\left( \ref{eq48}\right) $. Therefore
\begin{eqnarray*}
\int_{\Bbb{R}^N}\left| \Delta u\right| ^{q_1}dx &\leq &2^{q_1-1}\left[ \int_{%
\Bbb{R}^N}\left| V_1u\right| ^{q_1}dx+\left( \left\| V_2\right\| _{L^\infty
}+\left| \lambda \right| +C\right) ^{q_1}\int_{\Bbb{R}^N}\left| u\right|
^{q_1}dx\right] \\
&\leq &2^{q_1-1}\left( \int_{\Bbb{R}^N}\left| V_1\right| ^pdx\right) ^{\frac
1{\alpha _1}}\cdot \left( \int_{\Bbb{R}^N}\left| u\right| ^{q_1\beta
_1}dx\right) ^{\frac 1{\beta _1}} \\
&&+2^{q_1-1}\left( \left\| V_2\right\| _{L^\infty }+\left| \lambda \right|
+C\right) ^{q_1}\int_{\Bbb{R}^N}\left| u\right| ^{q_1}dx,
\end{eqnarray*}
$\alpha _1=\frac p{q_1}=1+\frac{\left( N-2q_0\right) p}{q_0N}$, $\beta _1=1+%
\frac{q_0N}{\left( N-2q_0\right) p}$. Notice that
\begin{equation}
q_1=\frac{q_0pN}{\left( N-2q_0\right) p+q_0N}\leq q_0^{*}\left( 2\right) =%
\frac{q_0N}{N-2q_0},  \label{eq49}
\end{equation}
and also
\begin{equation}
q_1\leq \frac{q_0pN}{\left( N-2q_0\right) p+q_0p}=\frac{q_0N}{N-q_0}
\label{eq50}
\end{equation}
via $p\leq N$, we derive
\begin{eqnarray}
\left\| \Delta u\right\| _{L^{q_1}} &\leq &2^{1-\frac 1{q_1}}\left\|
V_1\right\| _{L^p}\cdot \left\| u\right\| _{L^{q_0^{*}\left( 2\right)
}}+2^{1-\frac 1{q_1}}\left( \left\| V_2\right\| _{L^\infty }+\left| \lambda
\right| +C\right) \left\| u\right\| _{L^{q_1}}  \nonumber \\
\ &\leq &\theta \left( q_0,q_1,\lambda ,N\right) \cdot \left\| u\right\|
_{H^2},  \label{eq51}
\end{eqnarray}
and fully mimicking the trick employed by the preceding argument shows
\begin{equation}
\left\| u\right\| _{W^{2,q_1}}\leq \theta \left( q_0,q_1,\lambda
,N,\varepsilon \right) \left\| u\right\| _{H^2},  \label{eq52}
\end{equation}
where
\begin{eqnarray*}
\theta \left( q_0,q_1,\lambda ,N,\varepsilon \right) &:&=2\theta \left(
q_0,q_1,\lambda ,N\right) +2\left[ K\left( 2,q_1,N\right) \varepsilon
^{-1}+1\right] L\left( q_0,q_1,N\right) ; \\
\theta \left( q_0,q_1,\lambda ,N\right) &:&=2^{1-\frac 1{q_1}}\left[ \left\|
V_1\right\| _{L^p}C\left( 2,q_0,N\right) +\left( \left\| V_2\right\|
_{L^\infty }+\left| \lambda \right| +C\right) L\left( q_0,q_1,N\right)
\right] .
\end{eqnarray*}

Let $q_2=\frac{q_1pN}{\left( N-2q_1\right) p+q_1N}$, $\cdots $, $q_j=\frac{%
q_{j-1}pN}{\left( N-2q_{j-1}\right) p+q_{j-1}N}$. We presuppose $q_i<\frac
N2 $ for $0\leq i\leq j$, and so the validity of definitions is beyond doubt
on account of $\left( \ref{eq47}\right) $, and additionally $\left\{
q_i\right\} _{i=0}^j$ is monotonely increasing via $\left( \ref{eq48}\right)
$. Invoking backward induction we can get $q_i=\frac{2pN}{\left(
N-2iq_0\right) p+iq_0N}$. An argument analogous to case $\left\langle
a\right\rangle $ manifests that there is a $j_0\in \Bbb{N\cup }\left\{
0\right\} $, s.t. $q_{j_0}<\frac N2$, $q_{j_0+1}\geq \frac N2$, i.e. $\frac{%
p\left( N-4\right) }{4p-2N}-1\leq j_0<\frac{p\left( N-4\right) }{4p-2N}$.
Let
\begin{eqnarray*}
\Lambda ^{*} &=&\left\{ i\in \Bbb{N\cup }\left\{ 0\right\} :i<\frac{\left(
N-4\right) p}{4p-2N}\right\} ; \\
\Lambda _{*} &=&\left\{ i\in \Bbb{N\cup }\left\{ 0\right\} :i\geq \frac{%
\left( N-4\right) p}{4p-2N}-1\right\} .
\end{eqnarray*}
Obviously $\Lambda ^{*}\neq \varnothing $ and $\Lambda _{*}\neq \varnothing $%
. Observe that
\begin{equation}
\frac{\left( N-4\right) p}{4p-2N}>1\Leftrightarrow p\left( N-8\right) >-2N
\label{eq53}
\end{equation}
is established for $N\geq 6$, or $N=5$ with $p<\frac 23N$. Whereupon we
obtain that $j_0=l$ if $\frac{\left( N-4\right) p}{4p-2N}\in \left(
l,l+1\right) $, and $j_0=l-1$ if $\frac{\left( N-4\right) p}{4p-2N}=l$, and
here we denote $l=\left[ \frac{\left( N-4\right) p}{4p-2N}\right] $. Thus $%
j_0\geq 1$ if $N\geq 7$, or $p<N=6$, or $N=5$ with $p<\frac 23N$, and $j_0=0$
if $N=5$ with $p\geq \frac 23N$, or $p=N=6$.

In the case of $N\geq 7$, or $p<N=6$, or $N=5$ with $p<\frac 23N$, take
\[
\varepsilon <\Gamma \left( \frac 1{2K\left( 2,q_0,N\right) },\frac
1{2K\left( 2,q_1,N\right) },\cdots ,\frac 1{2K\left( 2,q_{j_0+1},N\right)
}\right) ,
\]
we yield
\begin{eqnarray}
\ \left\| u\right\| _{C_B^0} &=&\left\| u\right\| _{L^\infty }\leq C_0\left(
2,q_{j_0+1},N\right) \left\| u\right\| _{W^{2,q_{j_0+1}}}  \nonumber \\
&\leq &C_0\left( 2,q_{j_0+1},N\right) \theta \left(
q_{j_0},q_{j_0+1},\lambda ,N,\varepsilon \right) \left\| u\right\|
_{W^{2,q_{j_0}}}  \nonumber \\
&\leq &C_0\left( 2,q_{j_0+1},N\right) \prod\limits_{i=0}^{j_0}\theta \left(
q_i,q_{i+1},\lambda ,N,\varepsilon \right) \left\| u\right\| _{H^2}
\label{eq54}
\end{eqnarray}
followed by a routine argument provided $q_{j_0+1}>\frac N2$, where
\begin{eqnarray*}
\theta \left( q_i,q_{i+1},\lambda ,N,\varepsilon \right) &:&=2\theta \left(
q_i,q_{i+1},\lambda ,N\right) +2\left[ K\left( 2,q_{i+1},N\right)
\varepsilon ^{-1}+1\right] L\left( q_i,q_{i+1},N\right) ; \\
\theta \left( q_i,q_{i+1},\lambda ,N\right) &:&=2^{1-\frac 1{q_{i+1}}}\left[
\left\| V_1\right\| _{L^p}C\left( 2,q_i,N\right) +\left( \left\| V_2\right\|
_{L^\infty }+\left| \lambda \right| +C\right) L\left( q_i,q_{i+1},N\right)
\right] .
\end{eqnarray*}
And if $q_{j_0+1}=\frac N2$, observe that $W^{2,\frac N2}\left( \Bbb{R}%
^N\right) \hookrightarrow L^t\left( \Bbb{R}^N\right) $, $t\in \left[ \frac
N2,+\infty \right) $, we take $\widetilde{q}\in \left( q_{j_0+1},p\right) $,
hence
\begin{eqnarray*}
&&\left\| \Delta u\right\| _{L^{\widetilde{q}}} \\
&\leq &2^{1-\frac 1{^{\widetilde{q}}}}\left[ \left\| V_1\right\| _{L^p}\cdot
\left\| u\right\| _{L^{\frac{p\widetilde{q}}{p-\widetilde{q}}}}+\left(
\left\| V_2\right\| _{L^\infty }+\left| \lambda \right| +C\right) \left\|
u\right\| _{L^{\widetilde{q}}}\right] \\
&\leq &2^{1-\frac 1{^{\widetilde{q}}}}\left[ \left\| V_1\right\| _{L^p}%
\widetilde{C}\left( 2,q_{j_0+1},\frac{p\widetilde{q}}{p-\widetilde{q}}%
,N\right) +\left( \left\| V_2\right\| _{L^\infty }+\left| \lambda \right|
+C\right) \widetilde{C}\left( 2,q_{j_0+1},\widetilde{q},N\right) \right]
\cdot \left\| u\right\| _{W^{2,q_{j_0+1}}}.
\end{eqnarray*}
Let
\[
\varepsilon <\Gamma \left( \frac 1{2K\left( 2,q_0,N\right) },\frac
1{2K\left( 2,q_1,N\right) },\cdots ,\frac 1{2K\left( 2,q_{j_0+1},N\right)
},\frac 1{2K\left( 2,\widetilde{q},N\right) }\right) .
\]
Consequently,
\begin{eqnarray}
\ \left\| u\right\| _{C_B^0} &=&\left\| u\right\| _{L^\infty }\leq C_0\left(
2,\widetilde{q},N\right) \left\| u\right\| _{W^{2,\widetilde{q}}}  \nonumber
\\
&\leq &C_0\left( 2,\widetilde{q},N\right) \theta ^{*}\left( q_{j_0+1},%
\widetilde{q},\lambda ,N,\varepsilon \right) \left\| u\right\|
_{W^{2,q_{j_0+1}}}  \nonumber \\
&\leq &C_0\left( 2,\widetilde{q},N\right) \theta ^{*}\left( q_{j_0+1},p,%
\widetilde{q},\lambda ,N,\varepsilon \right) \prod\limits_{i=0}^{j_0}\theta
\left( q_i,q_{i+1},\lambda ,N,\varepsilon \right) \left\| u\right\| _{H^2},
\label{eq55}
\end{eqnarray}
where
\begin{eqnarray*}
\theta ^{*}\left( q_{j_0+1},p,\widetilde{q},\lambda ,N,\varepsilon \right)
&:&=2^{2-\frac 1{\widetilde{q}}}\left\| V_1\right\| _{L^p}\widetilde{C}%
\left( 2,q_{j_0+1},\frac{p\widetilde{q}}{p-\widetilde{q}},N\right) \\
&&+2^{2-\frac 1{\widetilde{q}}}\left( \left\| V_2\right\| _{L^\infty
}+\left| \lambda \right| +C\right) \widetilde{C}\left( 2,q_{j_0+1},%
\widetilde{q},N\right) \\
&&+2\left[ K\left( 2,\widetilde{q},N\right) \varepsilon ^{-1}+1\right]
\widetilde{C}\left( 2,q_{j_0+1},\widetilde{q},N\right) .
\end{eqnarray*}
As for the case $N=5$ with $p\geq \frac 23N$, or $p=N=6$, the consequence is
indeed more succinctly expressed by
\begin{equation}
\ \left\| u\right\| _{C_B^0}=\left\| u\right\| _{L^\infty }\leq C_0\left(
2,q_1,N\right) \left\| u\right\| _{W^{2,q_1}}\leq C_0\left( 2,q_1,N\right)
\theta \left( q_0,q_1,\lambda ,N,\varepsilon \right) \left\| u\right\| _{H^2}
\label{eq56}
\end{equation}
for $q_1>\frac N2$, and simultaneously
\begin{eqnarray}
\ \left\| u\right\| _{C_B^0} &=&\left\| u\right\| _{L^\infty }\leq C_0\left(
2,q,N\right) \left\| u\right\| _{W^{2,q}}  \nonumber \\
&\leq &C_0\left( 2,q,N\right) \theta ^{*}\left( q_1,p,q,\lambda
,N,\varepsilon \right) \theta \left( q_0,q_1,\lambda ,N,\varepsilon \right)
\left\| u\right\| _{H^2}.  \label{eq57}
\end{eqnarray}
for $q_1=\frac N2$ and $q\in \left( q_1,p\right) $. If $N=4$, the proof is
fairly trivial due to the emptiness of $\Lambda ^{*}$.

$\left\langle c\right\rangle $ $N\leq 3$. As $p=2$, we directly get $u\in
H^2\left( \Bbb{R}^N\right) \hookrightarrow C_B^1\left( \Bbb{R}^N\right) $ as
$N=1$, and $H^2\left( \Bbb{R}^N\right) \hookrightarrow L^q\left( \Bbb{R}%
^N\right) $ for $\forall q\in \left[ 2,\infty \right) $ as $N=2$, and $u\in
H^2\left( \Bbb{R}^N\right) \hookrightarrow C_B^0\left( \Bbb{R}^N\right) $ as
$N=3$. With this, the proof process terminates.\qed\vskip 5pt

We now go back to the question concerned. Based on Theorem\ref{thm2.2} we
can obtain the following consequence:

\begin{corollary}
\label{coroll2.3} Under the hypotheses of Proposition\ref{prop2.1}, assume $%
c\left( x\right) =c_1\left( x\right) +c_2\left( x\right) $, $c_1\left(
x\right) \in L^p\left( \Omega \right) $, $c_2\left( x\right) \in L^\infty
\left( \Omega \right) $, with $p=2$ if $N\leq 3$, $p>2$ as $N=4$, and $%
p>\frac N2$ as $N\geq 5$, $f\left( x\right) \in C\left( \overline{\Omega }%
\right) $. Then $\left( \ref{eq5}\right) $ exists a unique solution $u(x)\in
C^1\left( \overline{\Omega }\right) $ as $V_1=0$ a.e. on $\Omega $, and for
the case $\left\| V_1\right\| _{L^p}>0$, if $p>N\geq 4$ or $p=2$ and $N=1$, $%
u\in C^1\left( \overline{\Omega }\right) $, and if $\frac N2<p\leq N$ with $%
N\geq 4$ or $p=2$ and $N=2$,$3$, $u\in C\left( \overline{\Omega }\right) $.
Especially, if $c\left( x\right) $, $f\left( x\right) \in C^\alpha \left(
\overline{\Omega }\right) $, $0<\alpha <1$, then $u\left( x\right) \in
C^{2+\alpha }\left( \overline{\Omega }\right) $.
\end{corollary}

We omit the proof Corollary\ref{coroll2.3} considering that the process
carrys out substantially an argument analogous to Theorem\ref{thm2.2}.

\section{Compactness for the nonresonance case}

Aiming at verifying compactness, this section is mainly concerned with the
following nonlinear Schr\"odinger equation
\begin{equation}
\left\{
\begin{array}{cc}
Au=f\left( x,u\right) & x\in \Bbb{R}^N, \\
u\left( x\right) \rightarrow 0, & \text{as }\left| x\right| \rightarrow
\infty ,
\end{array}
\right.  \label{eq58}
\end{equation}
and the corresponding energetic functional of $\left( \ref{eq58}\right) $ is
of the form as:
\[
J\left( u\right) =J\left( u,a,b\right) =\frac 12\int_{\Bbb{R}^N}\left(
\left| \nabla u\right| ^2+V\left( x\right) u^2\right) dx-\int_{\Bbb{R}%
^N}F\left( x,u\right) \text{\textrm{, }}
\]
where $f\left( x,t\right) =au^{-}+bu^{+}+r\left( x,u\right) $, $F\left(
x,u\right) =\int_0^uf\left( x,s\right) ds$, $r\left( x,t\right) $ is a
Carath\'eodory function on $\Bbb{R}^N\times \Bbb{R}$, $\left( a,b\right)
\notin \Sigma \left( A\right) $, and $\Sigma \left( A\right) $ is called
Fu\v c\'\i k spectrum, defined as the set of all $\left( a,b\right) \in \Bbb{%
R}^2$ such that
\begin{equation}
\left\{
\begin{array}{cc}
Au=au^{-}+bu^{+}, & x\in \Bbb{R}^N, \\
u\left( x\right) \rightarrow 0, & \text{as\textrm{\ }}\left| x\right|
\rightarrow +\infty ,
\end{array}
\right.  \label{eq59}
\end{equation}
has a nontrivial solution $u$ in the form domain of $A$, where $A=-\Delta +V$%
, $u^{+}=\Lambda \left( u,0\right) $, $u^{-}=\Gamma \left( u,0\right) $.

Suppose:

$\left( r_1\right) $ $r\left( x,0\right) =0$, $\lim\limits_{\left| s\right|
\rightarrow \infty }\frac{r\left( x,s\right) }s=0$ uniformly with respect to
$x\in \Bbb{R}^N$.

$\left( r_2\right) $ $\exists \beta >0$, such that $\forall s_1$, $s_2\in
\Bbb{R}$, $\forall x\in \Bbb{R}^N$,
\[
\left| r\left( x,s_1\right) -r\left( x,s_2\right) \right| \leq \beta \left|
s_1-s_2\right|
\]
with $\beta <\sigma _0-\Lambda \left( a,b\right) $.

Under the assumptions $\left( r_1\right) \left( r_2\right) $, if $V$ is a
real K-R potential, $\sigma _{\text{\textrm{dis}}}\left( A\right) \neq
\varnothing $, $\inf \sigma \left( A\right) =\inf \sigma _{\text{\textrm{dis}%
}}\left( A\right) $, \cite{LiLi} proved that $J$ satisfies the (PS)
condition on $X$ for $\left( a,b\right) \notin \Sigma \left( A\right) $.
Actually, $\left( r_2\right) $ can be superseded by the hypothesis below:

$\left( \widetilde{r}_2\right) $ $\exists \beta >0$, such that $\forall s_1$%
, $s_2\in \Bbb{R}$, $\forall x\in \Bbb{R}^N$,
\[
\left| r\left( x,s_1\right) -r\left( x,s_2\right) \right| \leq \beta \left|
s_1-s_2\right| .
\]

We also need the additional hypotheses:

$\left( r_3\right) $ $\frac{r\left( x,s\right) }s\leq 0$ for $\forall s\in
\Bbb{R}\backslash \left\{ 0\right\} $, $\forall x\in \Bbb{R}^N$.

$\left( r_4\right) $ $r_0:=\lim\limits_{s\rightarrow 0}\frac{r\left(
x,s\right) }s$ uniformly with respect to $x\in \Bbb{R}^N$.

Our main theorem concerned compactness can be stated as follows:

\begin{theorem}
\label{thm3.1} Let $V$ be a real K-R potential, $\sigma _{\text{\textrm{dis}}%
}\left( A\right) \neq \varnothing $, $\inf \sigma \left( A\right) =\inf
\sigma _{\text{\textrm{dis}}}\left( A\right) $. Under the hypotheses $\left(
r_1\right) \left( \widetilde{r}_2\right) \left( r_3\right) \left( r_4\right)
$, if $\left( a,b\right) \notin \Sigma \left( A\right) $, $\Lambda \left(
a,b\right) <\sigma _0$, then $J$ satisfies the (PS) condition on $X$.
\end{theorem}

Let $h\left( x,s\right) :=r\left( x,s\right) s$ for $x\in \Bbb{R}^N$, $s\in
\Bbb{R}$. Expecting to verify Theorem\ref{thm3.1}, the following consequence
as a variant of Brezis-Lieb lemma is indispensable to back up the idea.

\begin{lemma}
\label{lem3.2} Let $\left\{ u_n\right\} _{n=1}^\infty \subset L^2\left( \Bbb{%
R}^N\right) $. Under the assumption $\left( \widetilde{r}_2\right) $, if $%
\left\{ u_n\right\} _{n=1}^\infty $ is bounded in $L^2\left( \Bbb{R}%
^N\right) $, $u_n\rightarrow u$ in $L_{\text{loc}}^2\left( \Bbb{R}^N\right) $%
, then
\begin{equation}
\lim\limits_{n\rightarrow \infty }\left( \int_{\Bbb{R}^N}h\left(
x,u_n\right) dx-\int_{\Bbb{R}^N}h\left( x,u_n-u\right) dx-\int_{\Bbb{R}%
^N}h\left( x,u\right) dx\right) =0\text{\textrm{.}}  \label{eq60}
\end{equation}
\end{lemma}

\textit{Proof. }Clearly $u\in L^2\left( \Bbb{R}^N\right) $ in terms of the
hypotheses. Making use of $\left( \widetilde{r}_2\right) $, we get
\begin{eqnarray}
&&h\left( x,u_n\right) -h\left( x,u_n-u\right) -h\left( x,u\right)  \nonumber
\\
&=&\left( r\left( x,u_n\right) -r\left( x,u_n-u\right) \right) \left(
u_n-u\right) ^{+}+\left( r\left( x,u_n\right) -r\left( x,u\right) \right)
u^{+}  \nonumber \\
&&+\left( r\left( x,u_n\right) -r\left( x,u_n-u\right) \right) \left(
u_n-u\right) ^{-}+\left( r\left( x,u_n\right) -r\left( x,u\right) \right)
u^{-}  \nonumber \\
&\geq &-\beta u^{+}\left( u_n-u\right) ^{+}+\beta u^{-}\left( u_n-u\right)
^{+}-\beta u^{+}\left( u_n-u\right) ^{+}+\beta u^{+}\left( u_n-u\right) ^{-}
\nonumber \\
&&+\beta u^{+}\left( u_n-u\right) ^{-}-\beta u^{-}\left( u_n-u\right)
^{-}+\beta \left( u_n-u\right) ^{+}u^{-}-\beta \left( u_n-u\right) ^{-}u^{-}
\nonumber \\
&=&-2\beta \left| u_n-u\right| \cdot \left| u\right|  \label{eq61}
\end{eqnarray}
and on the other side, an analogous argument shows
\begin{equation}
h\left( x,u_n\right) -h\left( x,u_n-u\right) -h\left( x,u\right) \leq 2\beta
\left| u_n-u\right| \cdot \left| u\right| \mathrm{.}  \label{eq62}
\end{equation}

Set $v_n=\left| u_n-u\right| $ and so $\left\{ v_n\right\} _{n=1}^\infty $
bounded in $L^2\left( \Bbb{R}^N\right) $. Assume $v_n\rightharpoonup u^{*}$
in $L^2\left( \Bbb{R}^N\right) $, and hence
\[
\int_\Omega v_n\varphi dx\rightarrow \int_\Omega u^{*}\varphi dx\text{, }%
\forall \varphi \in C_0^\infty \left( \Bbb{R}^N\right) .
\]
As $u_n\rightarrow u$ in $L_{\text{loc}}^2\left( \Bbb{R}^N\right) $,
\[
\int_\Omega v_n\varphi dx\rightarrow 0\text{, }\forall \varphi \in
C_0^\infty \left( \Bbb{R}^N\right) .
\]
Therefore, $u^{*}=0$ a.e on $\Bbb{R}^N$ and this alludes to
\begin{equation}
\int_{\Bbb{R}^N}v_n\cdot \left| u\right| dx\rightarrow 0.  \label{eq63}
\end{equation}
Combining $\left( \ref{eq61}\right) $, $\left( \ref{eq62}\right) $ and $%
\left( \ref{eq63}\right) $ we conclude the proof.\qed\vskip 5pt

\textit{Proof of Theorem3.1. }Due to\textit{\ }$\Lambda \left( a,b\right)
<\sigma _0$ and the hypothesis $\left( r_3\right) $, mimicking the trick
presented by \cite{LiLi} we can show that for any sequence $\left\{
u_k\right\} _{k=1}^\infty \subset X$, s.t. $\left\| J^{\prime }\left(
u_k\right) \right\| _X\rightarrow 0$, is bounded in $X$, so we concentrate
on substantiating the compactness, i.e., any sequence $\left\{ u_k\right\}
_{k=1}^\infty \subset X$, s.t. $\left\| J^{\prime }\left( u_k\right)
\right\| _X\rightarrow 0$, contains a convergent subsequence.

Since $\left\{ u_k\right\} _{k=1}^\infty $ is bounded in $X$, we assume $%
u_k\rightharpoonup u_0$ in $X$. For $\varphi \in C_0^\infty \left( \Bbb{R}%
^N\right) $, there is a bounded domain $\Omega \subset \Bbb{R}^N$, $\left\{
x\in \Bbb{R}^N:\varphi \left( x\right) \neq 0\right\} \subset \Omega $.
Hence, $u_k\rightarrow u_0$ in\textrm{\ }$L^2\left( \Omega \right) $. As
\begin{equation}
\left| f\left( x,s\right) \right| \leq \left( \Lambda \left( \left| a\right|
,\left| b\right| \right) +\beta \right) \cdot \left( 1+\left| s\right|
\right) ,  \label{eq64}
\end{equation}
for $\forall x\in \Bbb{R}^N$, $\forall s\in \Bbb{R}$, by TheoremA.2 \cite
{Willem},
\begin{eqnarray}
&&\ \ \left| \int_{\Bbb{R}^N}\left[ f\left( x,u_k\right) -f\left(
x,u_0\right) \right] \varphi dx\right|  \nonumber \\
\ &\leq &\left\| f\left( x,u_k\right) -f\left( x,u_0\right) \right\|
_{L^2\left( \Omega \right) }\cdot \left\| \varphi \right\| _{L^2\left(
\Omega \right) }\rightarrow 0,  \label{eq65}
\end{eqnarray}
and a standard argument shows
\begin{equation}
\lim\limits_{k\rightarrow \infty }\int_{\Bbb{R}^N}V\left( x\right)
u_k\varphi =\int_{\Bbb{R}^N}V\left( x\right) u_0\varphi .  \label{eq66}
\end{equation}
Consequently,
\begin{equation}
\int_{\Bbb{R}^N}\nabla u_0\nabla \varphi +V\left( x\right) u_0\varphi =\int_{%
\Bbb{R}^N}f\left( x,u_0\right) \varphi ,\text{\textrm{\ }}\forall \varphi
\in C_0^\infty \left( \Bbb{R}^N\right) .  \label{eq67}
\end{equation}

Based on the density of $C_0^\infty \left( \Bbb{R}^N\right) $ in $X$, we
yield
\begin{equation}
\int_{\Bbb{R}^N}\nabla u_0\nabla \varphi +V\left( x\right) u_0\varphi =\int_{%
\Bbb{R}^N}f\left( x,u_0\right) \varphi ,\text{\textrm{\ }}\forall \varphi
\in X,  \label{eq68}
\end{equation}
and this indicates that $u_0$ is a weak solution of $\left( \ref{eq58}%
\right) $.

Now we remain to prove $u_k\rightarrow u_0$ in $X$. Once again from
TheoremA.2 \cite{Willem},
\begin{equation}
r\left( x,u_k\right) \rightarrow r\left( x,u_0\right) \text{ in }L^2\left(
\Omega \right)  \label{eq69}
\end{equation}
for any bounded domain $\Omega \subset \Bbb{R}^N$, and thus
\begin{equation}
\int_{\Bbb{R}^N}r\left( x,u_k\right) \varphi \rightarrow \int_{\Bbb{R}%
^N}r\left( x,u_0\right) \varphi \text{\textrm{, }}\forall \varphi \in
C_0^\infty \left( \Bbb{R}^N\right) \text{\textrm{.}}  \label{eq70}
\end{equation}
Notice that $r\left( x,u_k\right) $ is bounded in $L^2\left( \Bbb{R}%
^N\right) $, so $r\left( x,u_k\right) \rightharpoonup v^{*}$ in $L^2\left(
\Bbb{R}^N\right) $. Hence
\begin{equation}
\int_{\Bbb{R}^N}r\left( x,u_k\right) \varphi \rightarrow \int_{\Bbb{R}%
^N}v^{*}\varphi \text{\textrm{, }}\forall \varphi \in C_0^\infty \left( \Bbb{%
R}^N\right) .  \label{eq71}
\end{equation}
$\left( \ref{eq70}\right) $ together with $\left( \ref{eq71}\right) $
derives $v^{*}=r\left( x,u_0\right) $. Thereby,
\begin{eqnarray}
&&\ \ \ \ \int_{\Bbb{R}^N}\left[ r\left( x,u_k\right) -r\left( x,u_0\right)
\right] \cdot \left( u_k-u_0\right) dx  \nonumber \\
\ &=&\int_{\Bbb{R}^N}h\left( x,u_k\right) -r\left( x,u_0\right) u_k-r\left(
x,u_k\right) u_0+h\left( x,u_0\right) dx  \nonumber \\
\ &=&\int_{\Bbb{R}^N}2h\left( x,u_0\right) +h\left( x,u_k-u_0\right)
-r\left( x,u_0\right) u_k-r\left( x,u_k\right) u_0dx+o\left( 1\right)
\nonumber \\
\ &=&\int_{\Bbb{R}^N}h\left( x,u_k-u_0\right) dx+o\left( 1\right) .
\label{eq72}
\end{eqnarray}
Being aware that
\begin{eqnarray}
&&\ \ \ \ \ \left\langle J^{\prime }\left( u_k\right) -J^{\prime }\left(
u_0\right) ,u_k-u_0\right\rangle _m  \nonumber \\
\ &=&\int_{\Bbb{R}^N}\left| \nabla \left( u_k-u_0\right) \right| ^2+V\left(
x\right) \left( u_k-u_0\right) ^2  \nonumber \\
&&\ \ \ \ \ \ \ \ \ \ -\int_{\Bbb{R}^N}\left[ f\left( x,u_k\right) -f\left(
x,u_0\right) \right] \cdot \left( u_k-u_0\right)  \nonumber \\
\ &\rightarrow &0  \label{eq73}
\end{eqnarray}
and inserting $\left( \ref{eq72}\right) $ into $\left( \ref{eq73}\right) $,
we obtain
\begin{eqnarray}
\left\| u_k-u_0\right\| _X^2 &\leq &\left( m+\Lambda \left( a,b\right)
\right) \left\| u_k-u_0\right\| _{L^2}^2+\int_{\Bbb{R}^N}h\left(
x,u_k-u_0\right) dx+o\left( 1\right)  \nonumber \\
\ &=&\left( m+\Lambda \left( a,b\right) \right) \left\| u_k-u_0\right\|
_{L^2}^2+\int_{\Bbb{R}^N}r\left( x,u_k-u_0\right) \left( u_k-u_0\right)
dx+o\left( 1\right)  \nonumber \\
\ &\leq &\left( m+\Lambda \left( a,b\right) \right) \left\| u_k-u_0\right\|
_{L^2}^2+o\left( 1\right)  \label{eq74}
\end{eqnarray}
via the hypotheses $\left( r_1\right) \left( r_3\right) $. Therefore we
deduce that $u_k\rightarrow u_0$ in $L^2\left( \Bbb{R}^N\right) $ and so $%
u_k\rightarrow u_0$ in $X$. The assertion follows.\qed\vskip 5pt

\section{A comparison theorem for the spectrum of Schr\"odinger operator}

Recall the comparison theorem for elliptic operators on a bounded domain $%
\Omega \subset \Bbb{R}^N$:

\begin{proposition}
\label{prop4.1} Let $\Omega \subset \Bbb{R}^N$ be a bounded domain with
smooth boundary, and $q_1\left( x\right) $, $q_2\left( x\right) \in C\left(
\overline{\Omega }\right) $, if $q_1\left( x\right) \leq q_2\left( x\right) $%
, then $\mu _j\left( -\Delta +q_1\right) \leq \mu _j\left( -\Delta
+q_2\right) $, and moreover if $q_1\left( x\right) \neq q_2\left( x\right) $%
, then $\mu _j\left( -\Delta +q_1\right) <\mu _j\left( -\Delta +q_2\right) $%
, where
\begin{eqnarray*}
\mu _j\left( -\Delta +q_1\right) :=\sup\limits_{\stackrel{\frak{B}%
_{j-1}\left( \Omega \right) \subset H_0^2\left( \Omega \right) }{\dim \frak{B%
}_{j-1}\left( \Omega \right) =j-1}}\inf\limits_{\stackrel{\phi \in \frak{B}%
_{j-1}^{\bot }\left( \Omega \right) \cap H_0^2\left( \Omega \right) }{%
\left\| \phi \right\| _{L^2\left( \Omega \right) }=1}}\left\langle A_1\phi
,\phi \right\rangle _{L^2}; \\
\mu _j\left( -\Delta +q_2\right) :=\sup\limits_{\stackrel{\frak{B}%
_{j-1}\left( \Omega \right) \subset H_0^2\left( \Omega \right) }{\dim \frak{B%
}_{j-1}\left( \Omega \right) =j-1}}\inf\limits_{\stackrel{\phi \in \frak{B}%
_{j-1}^{\bot }\left( \Omega \right) \cap H_0^2\left( \Omega \right) }{%
\left\| \phi \right\| _{L^2\left( \Omega \right) }=1}}\left\langle A_2\phi
,\phi \right\rangle _{L^2},
\end{eqnarray*}
where $\frak{B}_{j-1}^{\bot }\left( \Omega \right) $ is the orthogonal
complement in $L^2\left( \Omega \right) $ of $\frak{B}_{j-1}\left( \Omega
\right) $.
\end{proposition}

Clearly $\sigma \left( -\Delta +q_i\right) =\sigma _{\text{dis}}\left(
-\Delta +q_i\right) =\left\{ \mu _j\left( -\Delta +q_i\right) \right\}
_{j=1}^\infty $, $i=1$,$2$. A natural question arises: Does Proposition\ref
{prop4.1} hold for Schr\"odinger operator?\textbf{\ }A noticeable
counterexample is given by taking two K-R potentials $q_1$,$q_2$, s.t. $%
q_1>0 $ on $\Bbb{R}^N$ decaying to zero at infinity and $q_2=0$, since $%
\sigma \left( -\Delta +q_1\right) =\sigma \left( -\Delta +q_2\right) =\left[
0,\infty \right) $. Therefore, we need additional hypothesis aiming at
seeking for an analogue of Proposition\ref{prop4.1}. Our main consequence in
this section can be stated as follows:

\begin{theorem}
\label{thm4.2} Let $q_1$, $q_2$ be two K-R potentials, $q_1\leq q_2$, $%
q_1\neq q_2$. Let $A_i=-\Delta +q_i$, $i=1$,$2$. For fixed $j\in \Bbb{N}$,
define
\begin{eqnarray*}
\mu _j\left( A_1\right) :=\sup\limits_{\stackrel{\frak{B}_{j-1}\subset
H^2\left( \Bbb{R}^N\right) }{\dim \frak{B}_{j-1}=j-1}}\inf\limits_{\stackrel{%
\phi \in \frak{B}_{j-1}^{\bot }\cap H^2\left( \Bbb{R}^N\right) }{\left\|
\phi \right\| _{L^2}=1}}\left\langle A_1\phi ,\phi \right\rangle _{L^2}; \\
\mu _j\left( A_2\right) :=\sup\limits_{\stackrel{\frak{B}_{j-1}\subset
H^2\left( \Bbb{R}^N\right) }{\dim \frak{B}_{j-1}=j-1}}\inf\limits_{\stackrel{%
\phi \in \frak{B}_{j-1}^{\bot }\cap H^2\left( \Bbb{R}^N\right) }{\left\|
\phi \right\| _{L^2}=1}}\left\langle A_2\phi ,\phi \right\rangle _{L^2}.
\end{eqnarray*}
Set
\begin{eqnarray*}
\inf \sigma \left( A_1\right) =\mu _1\left( A_1\right) \leq \cdots \leq \mu
_k\left( A_1\right) \leq \inf \sigma _{\text{\textrm{ess}}}\left( A_1\right)
; \\
\inf \sigma \left( A_2\right) =\mu _1\left( A_2\right) \leq \cdots \leq \mu
_k\left( A_2\right) \leq \inf \sigma _{\text{\textrm{ess}}}\left( A_2\right)
,
\end{eqnarray*}
where $\frak{B}_{j-1}^{\bot }$ is the orthogonal complement in $L^2\left(
\Bbb{R}^N\right) $ of $\frak{B}_{j-1}$. Given $k\in \Bbb{N}$, if either
\[
\mu _k\left( A_1\right) \in \sigma \left( A_1\right) \backslash \sigma _{%
\text{\textrm{ess}}}\left( A_1\right) \text{ }
\]
or
\[
\mu _k\left( A_2\right) \in \sigma \left( A_2\right) \backslash \sigma _{%
\text{\textrm{ess}}}\left( A_2\right) ,\text{ }
\]
then
\begin{equation}
\mu _k\left( A_1\right) <\mu _k\left( A_2\right) .  \label{eq75}
\end{equation}
\end{theorem}

Before entering the proof of Theorem\ref{thm4.2}, we recall a min-max
principle of operator form:

\begin{proposition}
\label{prop4.3} (see TheoremXIII.1, \cite{RS})

Let $H$ be a self-adjoint operator on Hilbert space $E$ and $H$ is bounded
from below, i.e., $H\geq cI$ for some $c$. Define $\mu _k\left( H\right)
=\sup\limits_{\frak{B=}\text{\textrm{span}}\left\{ \psi _1,\cdots ,\psi
_{k-1}\right\} ,\psi _i\in D\left( H\right) }\inf\limits_{\stackrel{\psi \in
D\left( H\right) ,\left\| \psi \right\| _E=1}{\psi \in \frak{B}^{\bot }}%
}\left\langle H\psi ,\psi \right\rangle _E$. Then, for each fixed $k$,
either: $\left( a\right) $ there are $k$ eigenvalues (counting degenerate
eigenvalues a number of times equal to their multiplicity) below the bottom
of the essential spectrum, and $\mu _k\left( H\right) $ is the $k$-th
eigenvalue counting multiplicity; or $\left( b\right) $ $\mu _k$ is the
bottom of the essential spectrum, i.e., $\mu _k=\inf \left\{ \lambda
:\lambda \in \sigma _{\text{\textrm{ess}}}\left( H\right) \right\} $ and in
that case $\mu _k=\mu _{k+1}=\mu _{k+2}=\cdots $ and there are at most $k-1$
eigenvalues (counting multiplicity) below $\mu _k$.
\end{proposition}

\textit{Proof of Theorem4.2. }Take two cases into account:

$\left\langle i\right\rangle $ $\mu _k\left( A_1\right) \in \sigma \left(
A_1\right) \backslash \sigma _{\text{\textrm{ess}}}\left( A_1\right) $.

We need merely to deal with the case $\mu _k\left( A_2\right) \in \sigma _{%
\text{\textrm{ess}}}\left( A_2\right) $ and thus $\mu _k\left( A_2\right)
=\inf \sigma _{\text{\textrm{ess}}}\left( A_2\right) $. Suppose
\begin{eqnarray}
\inf \sigma _{\text{\textrm{ess}}}\left( A_1\right) &=&\mu _m\left(
A_1\right) :=\sup\limits_{\stackrel{\frak{B}_{j-1}\subset H^2\left( \Bbb{R}%
^N\right) }{\dim \frak{B}_{j-1}=j-1}}\inf\limits_{\stackrel{\phi \in \frak{B}%
_{j-1}^{\bot }\cap H^2\left( \Bbb{R}^N\right) }{\left\| \phi \right\|
_{L^2}=1}}\left\langle A_1\phi ,\phi \right\rangle _{L^2},  \nonumber \\
&=&\mu _{m+1}\left( A_1\right) =\mu _{m+2}\left( A_1\right) =\cdots ,
\label{eq76}
\end{eqnarray}
$m\in \Bbb{N}\cup \left\{ \infty \right\} $. Notice that
\begin{equation}
\mu _k\left( A_1\right) <\mu _m\left( A_1\right) \leq \mu _m\left(
A_2\right) \leq \inf \sigma _{\text{\textrm{ess}}}\left( A_2\right) =\mu
_k\left( A_2\right)  \label{eq77}
\end{equation}
by Proposition\ref{prop4.3}. The assertion follows.

$\left\langle ii\right\rangle $ $\mu _k\left( A_2\right) \in \sigma \left(
A_2\right) \backslash \sigma _{\text{\textrm{ess}}}\left( A_2\right) $.

By \cite{LiLi},
\begin{equation}
\mu _k\left( A_2\right) =\inf\limits_{\stackrel{\frak{B}_k\subset H^2}{\dim
\frak{B}_k=k}}\max\limits_{\stackrel{\psi \in \frak{B}_k}{\left\| \psi
\right\| _{L^2}=1}}\left\langle A_2\psi ,\psi \right\rangle _{L^2}.
\label{eq78}
\end{equation}
Let $\widetilde{E}_k=$\textrm{span}$\left\{ \widetilde{\psi }_1,\cdots ,%
\widetilde{\psi }_k\right\} $, where $\widetilde{\psi }_j$ is the
eigenvector corresponding to $\mu _j\left( A_2\right) $, $j=1$,$\cdots $,$k$%
. Thereby,
\begin{eqnarray}
\mu _k\left( A_2\right) &=&\max\limits_{\stackrel{\psi \in \widetilde{E}_k}{%
\left\| \psi \right\| _{L^2}=1}}\left\langle A_2\psi ,\psi \right\rangle
_{L^2}  \nonumber \\
&>&\max\limits_{\stackrel{\psi \in \widetilde{E}_k}{\left\| \psi \right\|
_{L^2}=1}}\left\langle A_1\psi ,\psi \right\rangle _{L^2}  \nonumber \\
&\geq &\inf\limits_{\stackrel{\frak{B}_k\subset H^2}{\dim \frak{B}_k=k}%
}\max\limits_{\stackrel{\psi \in \frak{B}_k}{\left\| \psi \right\| _{L^2}=1}%
}\left\langle A_1\psi ,\psi \right\rangle _{L^2}=\mu _k\left( A_1\right) .
\label{eq79}
\end{eqnarray}
The proof is complete.\qed\vskip 5pt

\section{Proof of Theorem1.1}

We divide the proof into four steps to secure the settled objective:

\textbf{Step 1. }We intend to produce a topological link with the
established assumptions. To see this, set $E_k=$\textrm{span}$\left\{ \psi
_1,\cdots ,\psi _k\right\} $, where $\psi _j$ is the eigenvector associated
with $\mu _j$, $j=1$,$\cdots $,$k$. We claim that for $\forall \varepsilon
>0 $, $\exists R>0$, for $u\in E_k$, $\left\| u\right\| _X\geq R$,
\begin{equation}
\int_{\Bbb{R}^N}\frac{\left| G\left( x,u\right) \right| }{\left\| u\right\|
_X^2}dx\leq \varepsilon .  \label{eq80}
\end{equation}
Argue indirectly, $\exists \varepsilon _0>0$, $u_n\in E_k$, $\left\|
u_n\right\| _X=\rho _n\rightarrow \infty $,
\begin{equation}
\lim\limits_{\overline{n\rightarrow \infty }}\int_{\Bbb{R}^N}\frac{\left|
G\left( x,u_n\right) \right| }{\left\| u_n\right\| _X^2}dx\geq \varepsilon
_0.  \label{eq81}
\end{equation}
Let $u_n=\rho _n\widehat{u}_n$, $\left\| \widehat{u}_n\right\| _X=1$ and we
can always choose the sequence $\left\{ \rho _n\right\} _{n=1}^\infty $ such
that $\left\{ \rho _n\right\} _{n=1}^\infty $ is strictly monotone
increasing with $\rho _n\rightarrow \infty $. Take $M_n=\rho _n^{\frac 12}$
and so $\left\{ \frac{M_n}{\rho _n}\right\} _{n=1}^\infty $ is strictly
monotone decreasing with $\frac{M_n}{\rho _n}\rightarrow 0$. Denote
\begin{eqnarray*}
\Theta _n^{\left( 1\right) } &:&=\left\{ x\in \Bbb{R}^N:\left| u_n\left(
x\right) \right| >M_n\right\} ; \\
\Theta _n^{\left( 2\right) } &:&=\left\{ x\in \Bbb{R}^N:\left| u_n\left(
x\right) \right| \leq M_n\right\} .
\end{eqnarray*}
Due to $\widehat{u}_n\rightarrow u^{*}$ in $E_k$, we define
\begin{eqnarray*}
\Gamma _n^{\left( 1\right) } &:&=\left\{ x\in \Bbb{R}^N:\left| \widehat{u}%
_n\left( x\right) \right| \leq \left| u^{*}\left( x\right) \right| \right\} ;
\\
\Gamma _n^{\left( 2\right) } &:&=\left\{ x\in \Bbb{R}^N:\left| \widehat{u}%
_n\left( x\right) \right| >\left| u^{*}\left( x\right) \right| \right\} .
\end{eqnarray*}
Notice that for $r_1$,$r_2\in \Bbb{R}$, $\left| r_1\right| \leq \left|
r_2\right| $,
\[
\left| G\left( x,r_1\right) \right| \leq \left| G\left( x,r_2\right) \right|
\]
as $0\leq r_1\leq r_2$ or $r_2\leq r_1\leq 0$, and
\[
\left| G\left( x,r_1\right) \right| \leq \left| G\left( x,-r_2\right)
\right|
\]
as $r_1\leq 0\leq r_2$ or $r_2\leq 0\leq r_1$. Let $r_1=\rho _n\widehat{u}%
_n\left( x\right) $, $r_2=\rho _nu^{*}\left( x\right) $, hence for $x\in
\Gamma _n^{\left( 1\right) }$,
\begin{equation}
\left| G\left( x,\rho _n\widehat{u}_n\right) \right| \leq \left| G\left(
x,\rho _nu^{*}\right) \right|  \label{eq82}
\end{equation}
as $0\leq \widehat{u}_n\leq u^{*}$ or $u^{*}\leq \widehat{u}_n\leq 0$, and
\begin{equation}
\left| G\left( x,\rho _n\widehat{u}_n\right) \right| \leq \left| G\left(
x,-\rho _nu^{*}\right) \right|  \label{eq83}
\end{equation}
as $\widehat{u}_n\leq 0\leq u^{*}$ or $u^{*}\leq 0\leq \widehat{u}_n$. Set
\begin{eqnarray*}
\Pi _n^{\left( 1\right) } &=&\left\{ x\in \Bbb{R}^N:0\leq \widehat{u}_n\leq
u^{*}\text{\textrm{\ or }}u^{*}\leq \widehat{u}_n\leq 0\right\} ; \\
\Pi _n^{\left( 2\right) } &=&\left\{ x\in \Bbb{R}^N:\widehat{u}_n\leq 0\leq
u^{*}\text{\textrm{\ or }}u^{*}\leq 0\leq \widehat{u}_n\right\} .
\end{eqnarray*}
Consequently,
\begin{eqnarray}
&&\ \int_{\Gamma _n^{\left( 1\right) }}\frac{\left| G\left( x,u_n\right)
\right| }{\rho _n^2}dx  \nonumber \\
\ &=&\int_{\Gamma _n^{\left( 1\right) }\cap \Pi _n^{\left( 1\right) }}\frac{%
\left| G\left( x,\rho _n\widehat{u}_n\right) \right| }{\rho _n^2}%
dx+\int_{\Gamma _n^{\left( 1\right) }\cap \Pi _n^{\left( 2\right) }}\frac{%
\left| G\left( x,\rho _n\widehat{u}_n\right) \right| }{\rho _n^2}dx
\nonumber \\
\ &\leq &\int_{\Gamma _n^{\left( 1\right) }\cap \Pi _n^{\left( 1\right) }}%
\frac{\left| G\left( x,\rho _nu^{*}\right) \right| }{\rho _n^2}%
dx+\int_{\Gamma _n^{\left( 1\right) }\cap \Pi _n^{\left( 2\right) }}\frac{%
\left| G\left( x,-\rho _nu^{*}\right) \right| }{\rho _n^2}dx  \label{eq84}
\end{eqnarray}
Let
\begin{eqnarray*}
\Sigma _n^{\left( 1\right) } &:&=\left\{ x\in \Bbb{R}^N:\left| \rho
_nu^{*}\right| >M_n\right\} ; \\
\Sigma _n^{\left( 2\right) } &:&=\left\{ x\in \Bbb{R}^N:\left| \rho
_nu^{*}\right| \leq M_n\right\} .
\end{eqnarray*}

In terms of the hypothesis $\left( g_7\right) $, $\forall \varepsilon >0$, $%
\exists N_1\in \Bbb{N}$, $\forall n\geq N_1$, $n\in \Bbb{N}$,
\begin{eqnarray}
\ \int_{\Gamma _n^{\left( 1\right) }\cap \Pi _n^{\left( 1\right) }\cap
\Sigma _n^{\left( 1\right) }}\frac{\left| G\left( x,\rho _nu^{*}\right)
\right| }{\rho _n^2}dx &\leq &\varepsilon \int_{\Bbb{R}^N}\left(
u^{*}\right) ^2dx,  \label{eq85} \\
\int_{\Gamma _n^{\left( 1\right) }\cap \Pi _n^{\left( 2\right) }\cap \Sigma
_n^{\left( 1\right) }}\frac{\left| G\left( x,-\rho _nu^{*}\right) \right| }{%
\rho _n^2}dx &\leq &\varepsilon \int_{\Bbb{R}^N}\left( u^{*}\right) ^2dx,
\label{eq86}
\end{eqnarray}
Also being aware that
\begin{eqnarray}
&&\ \ \ \ \ \ \ \ \int_{\Bbb{R}^N}\left( u^{*}\right) ^2dx  \nonumber \\
\ &=&\int_{\left| u^{*}\right| >\frac{M_1}{\rho _1}}\left( u^{*}\right)
^2dx+\sum\limits_{i=1}^\infty \int_{\frac{M_i}{\rho _i}\geq \left|
u^{*}\right| >\frac{M_{i+1}}{\rho _{i+1}}}\left( u^{*}\right) ^2dx  \nonumber
\\
\ &=&\int_{\left| u^{*}\right| >\frac{M_n}{\rho _n}}\left( u^{*}\right)
^2dx+\int_{\left| u^{*}\right| \leq \frac{M_n}{\rho _n}}\left( u^{*}\right)
^2dx,  \label{eq87}
\end{eqnarray}
so for $\forall \varepsilon >0$, $\exists N_2\in \Bbb{N}$, $\forall n\geq
N_2 $, $n\in \Bbb{N}$,
\begin{equation}
\int_{\left| u^{*}\right| \leq \frac{M_n}{\rho _n}}\left( u^{*}\right)
^2dx\leq \varepsilon ,  \label{eq88}
\end{equation}
and thus
\begin{eqnarray}
&&\ \ \ \ \ \ \ \int_{\Gamma _n^{\left( 1\right) }\cap \Pi _n^{\left(
1\right) }\cap \Sigma _n^{\left( 2\right) }}\frac{\left| G\left( x,\rho
_nu^{*}\right) \right| }{\rho _n^2}dx  \nonumber \\
\ &\leq &c\int_{\Gamma _n^{\left( 1\right) }\cap \Pi _n^{\left( 1\right)
}\cap \Sigma _n^{\left( 2\right) }}\frac{\left( \rho _nu^{*}\right) ^2}{\rho
_n^2}dx  \nonumber \\
\ &=&c\int_{\Gamma _n^{\left( 1\right) }\cap \Pi _n^{\left( 1\right) }\cap
\Sigma _n^{\left( 2\right) }}\left( u^{*}\right) ^2dx\leq c\varepsilon ,
\label{eq89}
\end{eqnarray}
and proceeding along the same lines,
\begin{equation}
\int_{\Gamma _n^{\left( 1\right) }\cap \Pi _n^{\left( 2\right) }\cap \Sigma
_n^{\left( 2\right) }}\frac{\left| G\left( x,-\rho _nu^{*}\right) \right| }{%
\rho _n^2}dx\leq c\varepsilon .  \label{eq90}
\end{equation}
Analogously, $\forall \varepsilon >0$, $\exists N_3\in \Bbb{N}$, for $n\geq
N_3$, $n\in \Bbb{N}$,
\begin{equation}
\int_{\Gamma _n^{\left( 2\right) }\cap \Theta _n^{\left( 1\right) }}\frac{%
\left| G\left( x,u_n\right) \right| }{\rho _n^2}dx\leq \varepsilon \int_{%
\Bbb{R}^N}\widehat{u}_n^2dx.  \label{eq91}
\end{equation}
A simple observation
\begin{eqnarray*}
&&\ \ \ \ \ \ \ \left| \int_{\Gamma _n^{\left( 2\right) }\cap \Theta
_n^{\left( 2\right) }}\widehat{u}_n^2-\left( u^{*}\right) ^2dx\right| \\
\ &\leq &\left\| \widehat{u}_n-u^{*}\right\| _{L^2}\cdot \left( \left\|
\widehat{u}_n\right\| _{L^2}+\left\| u^{*}\right\| _{L^2}\right) \rightarrow
0
\end{eqnarray*}
obtains that $\forall \varepsilon >0$, $\exists N_4\in \Bbb{N}$, for $n\geq
N_4$, $n\in \Bbb{N}$,
\begin{equation}
\int_{\Gamma _n^{\left( 2\right) }\cap \Theta _n^{\left( 2\right) }}\frac{%
\left| G\left( x,u_n\right) \right| }{\rho _n^2}dx\leq c\int_{\Gamma
_n^{\left( 2\right) }\cap \Theta _n^{\left( 2\right) }}\widehat{u}_n^2dx\leq
c\int_{\Gamma _n^{\left( 2\right) }\cap \Theta _n^{\left( 2\right) }}\left(
u^{*}\right) ^2dx+\varepsilon  \label{eq92}
\end{equation}
Take $N=\Lambda \left( N_1,N_2,N_3,N_4\right) $. Therefore, for $n\geq N$, $%
n\in \Bbb{N}$,
\begin{eqnarray}
&&\ \ \ \ \ \ \ \int_{\Bbb{R}^N}\frac{\left| G\left( x,u_n\right) \right| }{%
\left\| u_n\right\| _X^2}dx  \nonumber \\
\ &\leq &2\varepsilon \int_{\Bbb{R}^N}\left( u^{*}\right) ^2dx+2c\varepsilon
+\varepsilon \int_{\Bbb{R}^N}\widehat{u}_n^2dx+c\int_{\Gamma _n^{\left(
2\right) }\cap \Theta _n^{\left( 2\right) }}\left( u^{*}\right)
^2dx+\varepsilon ,  \label{eq93}
\end{eqnarray}
alluding to
\begin{equation}
\lim\limits_{n\rightarrow \infty }\int_{\Bbb{R}^N}\frac{\left| G\left(
x,u_n\right) \right| }{\left\| u_n\right\| _X^2}dx=0.  \label{eq94}
\end{equation}
No doubt we have made an erroneous assumption $\left( \ref{eq81}\right) $,
and the claim $\left( \ref{eq80}\right) $ is substantiated accordingly.

Let $\delta >0$, $\left[ \sigma _0-\delta ,\sigma _0\right) \cap \sigma
\left( A\right) =\varnothing $. In terms of $\left( \ref{eq80}\right) $, for
fixed $\widetilde{\varepsilon }\in \left( 0,\frac{\sigma _0-\delta -\mu _k}{%
\mu _k+m}\right) $, $\exists \widetilde{R}>0$, for $u\in E_k$, $\left\|
u\right\| _X\geq \widetilde{R}$,
\begin{eqnarray}
J_{\sigma _0-\delta }\left( u\right) &\leq &-\left( \frac{\sigma _0-\delta
-\mu _k}{\mu _k+m}-\widetilde{\varepsilon }\right) \left\| u\right\| _X^2
\nonumber \\
\ &\leq &-\left( \frac{\sigma _0-\delta -\mu _k}{\mu _k+m}-\widetilde{%
\varepsilon }\right) \widetilde{R}^2<0.  \label{eq95}
\end{eqnarray}
Notice that for the given $u\in X$, $J_\lambda \left( u\right) $ is
monotonely decreasing on $\lambda \in \Bbb{R}$, and therefore for $\lambda
\in \left( \sigma _0-\delta ,\sigma _0\right) $,
\begin{eqnarray}
\sup\limits_{u\in E_k\cap B\left( 0,\widetilde{R}\right) ^c}J_\lambda \left(
u\right) &\leq &\sup\limits_{u\in E_k\cap B\left( 0,\widetilde{R}\right)
^c}J_{\sigma _0-\delta }\left( u\right) \leq -\left( \frac{\sigma _0-\delta
-\mu _k}{\mu _k+m}-\widetilde{\varepsilon }\right) \widetilde{R}^2.
\nonumber \\
&&  \label{eq96}
\end{eqnarray}
On the other side, for $\lambda \in \left( \sigma _0-\delta ,\sigma
_0\right) $, $u\in E_k^{\bot }$,
\begin{equation}
J_\lambda \left( u\right) \geq -\int_{\Bbb{R}^N}G\left( x,u\right) dx\geq 0
\label{eq97}
\end{equation}
by $\left( g_3\right) $. Clearly $\partial D$ and $S$ homologically link(We
say that $\partial D$ and $S$ homologically link, if $\partial D\cap
S=\varnothing $ and $\left| \tau \right| \cap S\neq \varnothing $, for each
singular $k$ chain $\tau $ with $\widehat{\partial }\tau =\widehat{\partial }%
D$, where $\widehat{\partial }$ denotes the boundary operator and $\left|
\tau \right| $ is the support of $\tau $) by taking $S=E_k^{\bot }\cap X$, $%
D=B\left( 0,\widetilde{R}\right) \cap E_k$, where $E_k^{\bot }$ is the
orthogonal complement in $L^2\left( \Bbb{R}^N\right) $ of $E_k$.

It is easy to verify that $\exists \gamma ^{*}>0$, for $\forall \lambda \in
\left[ \sigma _0-\delta ,\sigma _0\right) $,
\begin{equation}
\sup\limits_{u\in \overline{D}}J_\lambda \left( u\right) \leq \gamma ^{*}
\label{eq98}
\end{equation}
by observing the two facts
\begin{equation}
\sup\limits_{u\in \overline{D}}J_\lambda \left( u\right) \leq
\sup\limits_{u\in \overline{D}}J_{\sigma _0-\delta }\left( u\right)
\label{eq99}
\end{equation}
and $J_{\sigma _0-\delta }$ is bounded in $\overline{D}$.

Set $\alpha ^{*}:=-\left( \frac{\sigma _0-\delta -\mu _k}{\mu _k+m}-%
\widetilde{\varepsilon }\right) \widetilde{R}^2$. Invoking Theorem1.1$%
^{\prime }$\cite{Cha}, for fixed $\beta _0>0$, $\forall \lambda \in \left[
\sigma _0-\delta ,\sigma _0\right) $,
\begin{equation}
H_k\left( J_\lambda ^{\gamma ^{*}+\beta _0},J_\lambda ^{\alpha ^{*}}\right)
\neq 0  \label{eq100}
\end{equation}
since $\partial D$ and $S$ homologically link, where $H_k\left( J_\lambda
^{\gamma ^{*}+\beta _0},J_\lambda ^{\alpha ^{*}}\right) $ denotes the
singular $k+1$-relative homology group of topological pair $\left( J_\lambda
^{\gamma ^{*}+\beta _0},J_\lambda ^{\alpha ^{*}}\right) $. By Theorem2.2\cite
{Cha}, there exists a critical point $u_\lambda $ with $M\left( J_\lambda
,u_\lambda \right) \geq k$. In view of the hypotheses of Theorem\ref{thm1.1}%
, $\lambda +g_0\leq \mu _i<\mu _k$ or $\lambda +g_0<\mu _i\leq \mu _k$ for $%
\lambda \in \left[ \sigma _0-\delta ,\sigma _0\right) $ as $\delta >0$ small
sufficiently. Set $V_\lambda \left( x\right) :=V\left( x\right) -\lambda $,
and thus
\begin{equation}
\mu _i\left( A_\lambda \right) :=\sup\limits_{\stackrel{\frak{B}%
_{i-1}\subset H^2}{\dim \frak{B}_{i-1}=i-1}}\inf\limits_{\stackrel{\phi \in
\frak{B}_{i-1}^{\bot }\cap H^2}{\left\| \phi \right\| _{L^2}=1}}\left\langle
A_\lambda \phi ,\phi \right\rangle _{L^2}>0,  \label{eq101}
\end{equation}
alluding to $M\left( J_\lambda ,0\right) \leq k-1$, where $A_\lambda
:=-\Delta +V_\lambda -g_0$. This yields $u_\lambda \neq 0$ and besides $%
c_\lambda :=J_\lambda \left( u_\lambda \right) >0$ by $\left( g_4\right) $.

\textbf{Step 2. }Arbitrarily choosing $\lambda _n\in \left( \sigma _0-\delta
,\sigma _0\right) $, and let $\lambda _n\rightarrow \sigma _0$ and denote by
$\left\{ u_{\lambda _n}\right\} _{n=1}^\infty $ the solution sequence
yielded by linking approach. In step2 we are devoted to showing that $%
\left\{ u_{\lambda _n}\right\} _{n=1}^\infty $ is bounded in $X$. We\textbf{%
\ }still use reduction to absurdity to derive the boundedness of $\left\{
u_{\lambda _n}\right\} _{n=1}^\infty $ in $X$. Suppose $\left\| u_{\lambda
_n}\right\| _X\rightarrow \infty $ and set $\widehat{u}_{\lambda _n}=\frac{%
u_{\lambda _n}}{\left\| u_{\lambda _n}\right\| _{H^1}}$. Then there exists a
renamed subsequence $\left\{ \widehat{u}_{\lambda _n}\right\} _{n=1}^\infty $
with $\widehat{u}_{\lambda _n}\rightharpoonup u_0$ solving
\begin{equation}
-\Delta \widehat{u}_{\lambda _n}+V\widehat{u}_{\lambda _n}=\lambda _n%
\widehat{u}_{\lambda _n}+\frac{g\left( x,u_{\lambda _n}\right) }{\left\|
u_{\lambda _n}\right\| _{H^1}}  \label{eq102}
\end{equation}
and satisfying the alternative:

1.(non-vanishing) there exist $\alpha >0$, $R>0$, s.t.
\begin{equation}
\lim\limits_{\overline{n\rightarrow +\infty }}\sup\limits_{y\in \Bbb{R}%
^N}\int_{B\left( y,R\right) }\left| \widehat{u}_{\lambda _n}\right| ^2dx\geq
\alpha >0.  \label{eq103}
\end{equation}

2.(vanishing) for $\forall R>0$%
\begin{equation}
\lim\limits_{n\rightarrow \infty }\sup\limits_{y\in \Bbb{R}^N}\int_{B\left(
y,R\right) }\left| \widehat{u}_{\lambda _n}\right| ^2dx=0.  \label{eq104}
\end{equation}

The non-vanishing case shows that $\exists N\in \Bbb{N}$, $\forall n\geq N$,
$\exists y_n\in \Bbb{R}^N$, $\left| y_n\right| \rightarrow \infty $,
\begin{equation}
\int_{B\left( y_n,R\right) }\left| \widehat{u}_{\lambda _n}\right| ^2dx\geq
\frac \alpha 2.  \label{eq105}
\end{equation}
Given $n\in \Bbb{N}$, set $x=z+y_n$ and the translated function $u^y\left(
z\right) =u\left( z+y\right) $. Hence
\[
u_{\lambda _n}\left( x\right) =u_{\lambda _n}\left( z+y_n\right) =u_{\lambda
_n}^{y_n}\left( z\right) .
\]
Define $v_n\left( z\right) :=u_{\lambda _n}^{y_n}\left( z\right) $, $%
\widehat{v}_n\left( z\right) =\frac{v_n\left( z\right) }{\left\| v_n\right\|
_{H^1}}$. Being aware that $\left\| v_n\right\| _{H^1}=\left\| u_{\lambda
_n}\right\| _{H^1}$, an alternative to $\left( \ref{eq105}\right) $ can be
easily seen as
\begin{equation}
\int_{B\left( 0,R\right) }\left| \widehat{v}_n\left( z\right) \right|
^2dz\geq \frac \alpha 2,  \label{eq106}
\end{equation}
so $\widehat{v}_n\rightharpoonup v_0\neq 0$ in $H^1\left( \Bbb{R}^N\right) $
by $\left( \ref{eq106}\right) $.

The equation
\begin{equation}
-\Delta \widehat{u}_{\lambda _n}\left( x\right) +V\left( x\right) \widehat{u}%
_{\lambda _n}\left( x\right) =\lambda _n\widehat{u}_{\lambda _n}\left(
x\right) +\frac{g\left( x,u_{\lambda _n}\left( x\right) \right) }{\left\|
u_{\lambda _n}\right\| _{H^1}}  \label{eq107}
\end{equation}
can be rewritten as
\begin{equation}
-\Delta \widehat{v}_n\left( z\right) +V\left( z+y_n\right) \widehat{v}%
_n\left( z\right) =\lambda _n\widehat{v}_n\left( z\right) +\frac{g\left(
z+y_n,v_n\left( z\right) \right) }{\left\| v_n\right\| _{H^1}}  \label{eq108}
\end{equation}

Clearly, for $\forall \varphi \in C_0^\infty \left( \Bbb{R}^N\right) $,
\begin{equation}
\lim\limits_{n\rightarrow \infty }\frac{\int_{\Bbb{R}^N}g\left(
z+y_n,v_n\left( z\right) \right) \varphi \left( z\right) dz}{\left\|
v_n\right\| _{H^1}}=0  \label{eq109}
\end{equation}
by \cite{LiLi}. And also notice that
\begin{eqnarray}
&&\ \ \int_{\Bbb{R}^N}V\left( z+y_n\right) \widehat{v}_n\left( z\right)
\varphi \left( z\right) dz  \nonumber \\
\ &=&\int_\Omega V\left( z+y_n\right) \widehat{v}_n\left( z\right) \varphi
\left( z\right) dz  \nonumber \\
\ &\rightarrow &\sigma _0\int_\Omega v_0\left( z\right) \varphi \left(
z\right) dz=\sigma _0\int_{\Bbb{R}^N}v_0\left( z\right) \varphi \left(
z\right) dz  \label{eq110}
\end{eqnarray}
by \cite{JT}, combining $\left( \ref{eq108}\right) $, $\left( \ref{eq109}%
\right) $ and $\left( \ref{eq110}\right) $, we have
\begin{equation}
-\Delta v_0\left( z\right) =0.  \label{eq111}
\end{equation}
This manifests that $0$ is an eigenvalue of $-\Delta $ with the eigenvector $%
v_0\left( z\right) $ and thereby facilitates a paradox stemming from $\sigma
\left( -\Delta \right) =\sigma _c\left( -\Delta \right) =\left[ 0,\infty
\right) $, where $\sigma _c\left( -\Delta \right) $ stands for continuous
spectrum of $-\Delta $. The non-vanishing case is impossible.

Given $M>0$, define $\Omega _{n,M}:=\left\{ x\in \Bbb{R}^N:\left| u_{\lambda
_n}\left( x\right) \right| \geq M\right\} $. We claim $\left| \Omega
_{n,M}\right| \rightarrow \infty $ as $n\rightarrow \infty $. By way of
negation, there exists $\gamma _0>0$ with $\left| \Omega _{n,M}\right| \leq
\gamma _0$ for $n\in \Bbb{N}$. It is pretty straightforward to check that
\begin{eqnarray}
\ \ \left\langle A_m\widehat{u}_{\lambda _n},\widehat{u}_{\lambda
_n}\right\rangle _{L^2} &=&\left( \lambda _n+m\right) \int_{\Bbb{R}^N}%
\widehat{u}_{\lambda _n}^2dx+\int_{\Bbb{R}^N}\frac{g\left( x,u_{\lambda
_n}\right) }{u_{\lambda _n}}\widehat{u}_{\lambda _n}^2dx  \nonumber \\
\ &\leq &\left( \lambda _n+m\right) \int_{\Bbb{R}^N}\widehat{u}_{\lambda
_n}^2dx+\int_{\Bbb{R}^N\backslash \Omega _{n,M}}\frac{g\left( x,u_{\lambda
_n}\right) }{u_{\lambda _n}}\widehat{u}_{\lambda _n}^2dx  \nonumber \\
\ &\leq &\left( \lambda _n+m\right) \int_{\Bbb{R}^N}\widehat{u}_{\lambda
_n}^2dx-C_M\int_{\Bbb{R}^N}\widehat{u}_{\lambda _n}^2dx+C_M\int_{\Omega
_{n,M}}\widehat{u}_{\lambda _n}^2dx  \nonumber \\
\ &=&\left( \lambda _n+m-C_M\right) \int_{\Bbb{R}^N}\widehat{u}_{\lambda
_n}^2dx+C_M\int_{\Omega _{n,M}}\widehat{u}_{\lambda _n}^2dx  \label{eq112}
\end{eqnarray}
since we can easily infer from $\left( g_3\right) $ that there exists $C_M>0$%
, s.t., $g_0\leq \frac{g\left( x,t\right) }t\leq -C_M$ for $x\in \Bbb{R}^N$,
$t\in \left[ -M,M\right] $. Hence $\left( \ref{eq112}\right) $ yields
\begin{eqnarray}
&&\ \ \left( \mu _1+m\right) \left\| P_{E_k}\widehat{u}_{\lambda _n}\right\|
_{L^2}^2+\left( \sigma _0+m\right) \left\| P_{E_k^{\bot }}\widehat{u}%
_{\lambda _n}\right\| _{L^2}^2  \nonumber \\
\ &\leq &\int_{\Bbb{R}^N}\left| \nabla P_{E_k}\widehat{u}_{\lambda
_n}\right| ^2+\left( V\left( x\right) +m\right) \left( P_{E_k}\widehat{u}%
_{\lambda _n}\right) ^2  \nonumber \\
&&\ \ +\int_{\Bbb{R}^N}\left| \nabla P_{E_k^{\bot }}\widehat{u}_{\lambda
_n}\right| ^2+\left( V\left( x\right) +m\right) \left( P_{E_k^{\bot }}%
\widehat{u}_{\lambda _n}\right) ^2  \nonumber \\
\ &\leq &\left( \lambda _n+m-C_M\right) \int_{\Bbb{R}^N}\widehat{u}_{\lambda
_n}^2dx+C_M\left( \int_{\Bbb{R}^N}\left| \widehat{u}_{\lambda _n}\right|
^pdx\right) ^{\frac 2p}\gamma _0^{1-\frac 2p}  \nonumber \\
\ &=&\left( \lambda _n+m-C_M\right) \left\| P_{E_k}\widehat{u}_{\lambda
_n}\right\| _{L^2}^2+\left( \lambda _n+m-C_M\right) \left\| P_{E_k^{\bot }}%
\widehat{u}_{\lambda _n}\right\| _{L^2}^2  \nonumber \\
&&\ \ +C_M\left( \int_{\Bbb{R}^N}\left| \widehat{u}_{\lambda _n}\right|
^pdx\right) ^{\frac 2p}\gamma _0^{1-\frac 2p}  \label{eq113}
\end{eqnarray}
for $p\in \left( 2,2^{*}\right) $ via H\"older inequality, and consequently
derives
\begin{eqnarray}
&&\ \ \left\| P_{E_k^{\bot }}\widehat{u}_{\lambda _n}\right\| _{L^2}^2
\nonumber \\
\ &\leq &\frac{\lambda _n-\mu _1-C_M}{\sigma _0-\lambda _n+C_M}\left\|
P_{E_k}\widehat{u}_{\lambda _n}\right\| _{L^2}^2+\frac{C_M}{\sigma
_0-\lambda _n+C_M}\left( \int_{\Bbb{R}^N}\left| \widehat{u}_{\lambda
_n}\right| ^pdx\right) ^{\frac 2p}\gamma _0^{1-\frac 2p}  \label{eq114}
\end{eqnarray}
for $n\in \Bbb{N}$ sufficiently large, where $P_{E_k}$ is the orthogonal
projection onto $E_k$ in $L^2\left( \Bbb{R}^N\right) $. Due to $\widehat{u}%
_{\lambda _n}\rightharpoonup 0$ in $X$, $P_{E_k}\widehat{u}_{\lambda
_n}\rightarrow 0$ in $X$ and thus $P_{E_k}\widehat{u}_{\lambda
_n}\rightarrow 0$ in $L^2\left( \Bbb{R}^N\right) $. By \cite{Lions}, the
vanishing case obtains directly that $\left\| \widehat{u}_{\lambda
_n}\right\| _{L^r}\rightarrow 0$ for $\forall r\in \left( 2,2^{*}\right) $.
Accordingly,
\begin{equation}
\left\| P_{E_k^{\bot }}\widehat{u}_{\lambda _n}\right\| _{L^2}^2\rightarrow 0
\label{eq115}
\end{equation}
and so we get $\left\| \widehat{u}_{\lambda _n}\right\| _{L^2}\rightarrow 0$%
, which gives rise to a self-contradictory inequality
\begin{equation}
\ 1=\left\langle A_m\widehat{u}_{\lambda _n},\widehat{u}_{\lambda
_n}\right\rangle _{L^2}=\left( \lambda _n+m\right) \int_{\Bbb{R}^N}\widehat{u%
}_{\lambda _n}^2dx+\int_{\Bbb{R}^N}\frac{g\left( x,u_{\lambda _n}\right) }{%
u_{\lambda _n}}\widehat{u}_{\lambda _n}^2dx\rightarrow 0.  \label{eq116}
\end{equation}
We conclude the claim as predicted.

One deduces from $\left( g_4\right) $ that there exists a $b_M>0$,
\begin{equation}
g\left( x,t\right) t-2G\left( x,t\right) >b_M  \label{eq117}
\end{equation}
for $x\in \Bbb{R}^N$, $\left| t\right| \geq M$. Therefore, an absurd
assertion occurs as
\begin{eqnarray}
c &\geq &J_{\lambda _n}\left( u_{\lambda _n}\right)  \nonumber \\
\ &=&\frac 12\left\| u_{\lambda _n}\right\| _X^2-\frac{\lambda _n+m}2\int_{%
\Bbb{R}^N}u_{\lambda _n}^2dx-\int_{\Bbb{R}^N}G\left( x,u_{\lambda _n}\right)
dx  \nonumber \\
\ &=&\int_{\Bbb{R}^N}\left[ \frac 12g\left( x,u_{\lambda _n}\right)
u_{\lambda _n}-G\left( x,u_{\lambda _n}\right) \right] dx  \nonumber \\
\ &\geq &\frac{b_M}2\left| \Omega _{n,M}\right| \rightarrow \infty .
\label{eq118}
\end{eqnarray}
Thus far we have validated the prophecy that $\left\{ u_{\lambda _n}\right\}
_{n=1}^\infty $ is bounded in $X$.

\textbf{Step 3}. Suppose $u_{\lambda _n}\rightharpoonup u^{*}$ in $X$ and so
$u^{*}\in K_{J_{\sigma _0}}$. We are left to substantiate that $u_{\lambda
_n}\rightarrow u^{*}$ in $X$. Indeed Theorem\ref{thm2.2} directly exhibits
the boundedness of $\left\{ u_{\lambda _n}\right\} _{n=1}^\infty $ in $%
L^\infty \left( \Bbb{R}^N\right) $ (The proof for the cases $\frac N2<p\leq
N $ with $N\geq 4$, and $N\leq 3$ will be put into the appendix). As a
consequence, there exists a $\beta >0$,
\begin{equation}
\Lambda \left( \sup\limits_{n\in \Bbb{N}}\left\| u_{\lambda _n}\right\|
_{L^\infty },\left\| u^{*}\right\| _{L^\infty }\right) \leq \beta .
\label{eq119}
\end{equation}
Once again from Lemma\ref{lem3.2}, there exists a $\rho _0>0$
\begin{eqnarray}
&&\ \ \ \left\langle A_m\left( u_{\lambda _n}-u^{*}\right) ,u_{\lambda
_n}-u^{*}\right\rangle _{L^2}  \nonumber \\
\ &=&\left( \sigma _0+m\right) \left\| u_{\lambda _n}-u^{*}\right\|
_{L^2}^2+\int_{\Bbb{R}^N}\left[ g\left( x,u_{\lambda _n}\right) -g\left(
x,u^{*}\right) \right] \left( u_{\lambda _n}-u^{*}\right) dx  \nonumber \\
&&\ \ \ +\left( \lambda _n-\sigma _0\right) \int_{\Bbb{R}^N}u_{\lambda
_n}\left( u_{\lambda _n}-u^{*}\right) dx  \nonumber \\
\ &=&\left( \sigma _0+m\right) \left\| u_{\lambda _n}-u^{*}\right\|
_{L^2}^2+\int_{\Bbb{R}^N}\frac{g\left( x,u_{\lambda _n}-u^{*}\right) }{%
u_{\lambda _n}-u^{*}}\left( u_{\lambda _n}-u^{*}\right) ^2dx  \nonumber \\
&&\ \ +\left( \lambda _n-\sigma _0\right) \int_{\Bbb{R}^N}u_{\lambda
_n}\left( u_{\lambda _n}-u^{*}\right) dx+o\left( 1\right)  \nonumber \\
\ &\leq &\left( \sigma _0-\rho _0+m\right) \left\| u_{\lambda
_n}-u^{*}\right\| _{L^2}^2+\left( \sigma _0-\lambda _n\right) \left\|
u_{\lambda _n}\right\| _{L^2}\cdot \left\| u_{\lambda _n}-u^{*}\right\|
_{L^2}+o\left( 1\right) .  \nonumber \\
&&\ \ \   \label{eq120}
\end{eqnarray}
A standard argument on $\left( \ref{eq120}\right) $ produces
\begin{eqnarray}
&&\ \ \ \left\| P_{E_k^{\bot }}\left( u_{\lambda _n}-u^{*}\right) \right\|
_{L^2}^2  \nonumber \\
\ &\leq &\frac{\sigma _0-\mu _1-\rho _0}{\rho _0}\left\| P_{E_k}\left(
u_{\lambda _n}-u^{*}\right) \right\| _{L^2}^2+\frac{\sigma _0-\lambda _n}{%
\rho _0}\left\| u_{\lambda _n}\right\| _{L^2}\cdot \left\| u_{\lambda
_n}-u^{*}\right\| _{L^2}+o\left( 1\right)  \nonumber \\
\ &\rightarrow &0,  \label{eq121}
\end{eqnarray}
alluding to $u_{\lambda _n}\rightarrow u^{*}$ in $X$. The compactness is
corroborated accordingly.

\textbf{Step 4}. Thanks to $\sigma _0+g_0\leq \mu _i<\mu _k$ or $\sigma
_0+g_0<\mu _i$, we get $M\left( J_{\sigma _0},0\right) \leq k-1$. Also
observe that $M\left( J_{\lambda _n},u_{\lambda _n}\right) \geq k$ for fixed
$n\in \Bbb{N}$, and $u_{\lambda _n}\rightarrow u^{*}$ in $X$, so $M\left(
J_{\sigma _0},u^{*}\right) \geq k$. Therefore $u^{*}\neq 0$.

Next we are devoted to verifying that if $g_0<0$ and $\left( g_4\right) $
holds, then an aibitrary nonzero solution $u_\lambda $ of $J_\lambda $ is
nondegenerate with
\begin{equation}
m\left( J_\lambda ,u_\lambda \right) =M\left( J_\lambda ,u_\lambda \right) =k
\label{eq122}
\end{equation}
for the case $K_{J_\lambda }\backslash \left\{ 0\right\} \neq \varnothing $,
$\mu _{k-1}<\lambda +g_0<\mu _k<\lambda \leq \sigma _0$.

Given $\lambda \in \Bbb{R}$, and
\[
u\in D\left( \Bbb{R}^N\right) :=\left\{ u\in X:\text{mes}\left\{ x\in \Bbb{R}%
^N:u\left( x\right) =0\right\} =0\right\} ,
\]
we set $\widetilde{A}_{\lambda ,u}=-\Delta +\widetilde{V}_{\lambda ,u}$ with
\[
\widetilde{V}_{\lambda ,u}\left( x\right) :=\left\{
\begin{array}{cc}
V_\lambda \left( x\right) -g_0, & x\in \Bbb{R}^N,\text{\textrm{\ }}u\left(
x\right) =0, \\
V_\lambda \left( x\right) -\frac{g\left( x,u\left( x\right) \right) }{%
u\left( x\right) }, & x\in \Bbb{R}^N,\text{\textrm{\ }}u\left( x\right) \neq
0,
\end{array}
\right.
\]
where $V_\lambda \left( x\right) =V\left( x\right) -\lambda $. Let $%
u_\lambda $ be the nonzero solution of $\left( \ref{eq1}\right) $. As $%
V_\lambda \in L_{loc}^{\frac N2}\left( \Bbb{R}^N\right) $, by the unique
continuation property of solutions of Schr\"odinger equation(see \cite{Heinz}%
, \cite{JK}, \cite{K}, \cite{Niren2}), we get $u_\lambda \in D\left( \Bbb{R}%
^N\right) $. We now claim that $0$ is the $j$-th eigenvalue of $\widetilde{A}%
_{\lambda ,u_\lambda }=-\Delta +\widetilde{V}_{\lambda ,u_\lambda }$ with $%
j=k$, where
\[
\widetilde{V}_{\lambda ,u_\lambda }\left( x\right) :=\left\{
\begin{array}{cc}
V_\lambda \left( x\right) -g_0, & x\in \Bbb{R}^N,\text{\textrm{\ }}u\left(
x\right) =0, \\
V_\lambda \left( x\right) -\frac{g\left( x,u_\lambda \left( x\right) \right)
}{u_\lambda \left( x\right) }, & x\in \Bbb{R}^N,\text{\textrm{\ }}u\left(
x\right) \neq 0.
\end{array}
\right.
\]
Divide the proof into two cases:

$\left\langle i\right\rangle $ $\lambda <\sigma _0$.

Indeed, due to
\begin{eqnarray}
&&\ \ \ \ \ \ \sup\limits_{\stackrel{\frak{B}_{k-2}\subset H^2}{\dim \frak{B}%
_{k-2}=k-2}}\inf\limits_{\stackrel{\phi \in \frak{B}_{k-2}^{\bot }\cap H^2}{%
\left\| \phi \right\| _{L^2}=1}}\left\langle \widetilde{A}_{\lambda
,u_\lambda }\phi ,\phi \right\rangle _{L^2}  \nonumber \\
\ &\leq &\sup\limits_{\stackrel{\frak{B}_{k-2}\subset H^2}{\dim \frak{B}%
_{k-2}=k-2}}\inf\limits_{\stackrel{\phi \in \frak{B}_{k-2}^{\bot }\cap H^2}{%
\left\| \phi \right\| _{L^2}=1}}\left\langle A_\lambda \phi ,\phi
\right\rangle _{L^2}=\mu _{k-1}-\lambda -g_0<0,  \label{eq123}
\end{eqnarray}
we obtain $\mu _{k-1}\left( \widetilde{A}_{\lambda ,u_\lambda }\right) <0$,
where $A_\lambda =-\Delta +V_\lambda -g_0$. On the other side, notice that
\begin{equation}
-\Delta +V_\lambda \leq \widetilde{A}_{\lambda ,u_\lambda },  \label{eq124}
\end{equation}
hence
\begin{eqnarray}
\mu _{k+1}\left( \widetilde{A}_{\lambda ,u_\lambda }\right) &=&\sup\limits_{%
\stackrel{\frak{B}_k\subset H^2}{\dim \frak{B}_k=k}}\inf\limits_{\stackrel{%
\phi \in \frak{B}_k^{\bot }\cap H^2}{\left\| \phi \right\| _{L^2}=1}%
}\left\langle \widetilde{A}_{\lambda ,u_\lambda }\phi ,\phi \right\rangle
_{L^2}  \nonumber \\
&\geq &\sup\limits_{\stackrel{\frak{B}_k\subset H^2}{\dim \frak{B}_k=k}%
}\inf\limits_{\stackrel{\phi \in \frak{B}_k^{\bot }\cap H^2}{\left\| \phi
\right\| _{L^2}=1}}\left\langle \left( -\Delta +V_\lambda \right) \phi ,\phi
\right\rangle _{L^2}  \nonumber \\
&=&\sigma _0-\lambda .  \label{eq125}
\end{eqnarray}
The assertion follows.

$\left\langle ii\right\rangle $ $\lambda =\sigma _0$.

If
\begin{equation}
\sup\limits_{\stackrel{\frak{B}_k\subset H^2}{\dim \frak{B}_k=k}%
}\inf\limits_{\stackrel{\phi \in \frak{B}_k^{\bot }\cap H^2}{\left\| \phi
\right\| _{L^2}=1}}\left\langle \widetilde{A}_{\lambda ,u_\lambda }\phi
,\phi \right\rangle _{L^2}<\inf \sigma _{\text{ess}}\left( \widetilde{A}%
_{\lambda ,u_\lambda }\right) =\sigma _0-\lambda -g_0,  \label{eq126}
\end{equation}
Making use of Theorem\ref{thm4.2}, the first `` $\geq $ '' in $\left( \ref
{eq125}\right) $ can be substituted by `` $>$'', and therefore
\begin{equation}
\mu _{k+1}\left( \widetilde{A}_{\lambda ,u_\lambda }\right) =\sup\limits_{%
\stackrel{\frak{B}_k\subset H^2}{\dim \frak{B}_k=k}}\inf\limits_{\stackrel{%
\phi \in \frak{B}_k^{\bot }\cap H^2}{\left\| \phi \right\| _{L^2}=1}%
}\left\langle \widetilde{A}_{\lambda ,u_\lambda }\phi ,\phi \right\rangle
_{L^2}>0.  \label{eq127}
\end{equation}
Alternatively if
\begin{equation}
\sup\limits_{\stackrel{\frak{B}_k\subset H^2}{\dim \frak{B}_k=k}%
}\inf\limits_{\stackrel{\phi \in \frak{B}_k^{\bot }\cap H^2}{\left\| \phi
\right\| _{L^2}=1}}\left\langle \widetilde{A}_{\lambda ,u_\lambda }\phi
,\phi \right\rangle _{L^2}=\sigma _0-\lambda -g_0,  \label{eq128}
\end{equation}
we immediately derive $\mu _{k+1}\left( \widetilde{A}_{\lambda ,u_\lambda
}\right) >0$.

Since $0\in \sigma _{\text{dis}}\left( \widetilde{A}_{\lambda ,u_\lambda
}\right) \backslash \sigma _{\text{ess}}\left( \widetilde{A}_{\lambda
,u_\lambda }\right) $, the claim follows. Let $A_{\lambda ,u_\lambda
}:=-\Delta +V_\lambda -g_t^{\prime }\left( x,u_\lambda \right) $. Notice
that $\left( g_4\right) \Rightarrow A_{\lambda ,u_\lambda }\leq \widetilde{A}%
_{\lambda ,u_\lambda }$, we yield
\begin{eqnarray}
\mu _k\left( \widetilde{A}_{\lambda ,u_\lambda }\right) &\leq &\mu _k\left(
A_\lambda \right) \leq \mu _k\left( A\right) -\lambda -g_0  \nonumber \\
&<&\sigma _0-\lambda -g_0=\inf \sigma _{\text{ess}}\left( \widetilde{A}%
_{\lambda ,u_\lambda }\right)  \label{eq129}
\end{eqnarray}
and consequently
\begin{equation}
\mu _k\left( A_{\lambda ,u_\lambda }\right) <\mu _k\left( \widetilde{A}%
_{\lambda ,u_\lambda }\right) =0  \label{eq130}
\end{equation}
via Theorem\ref{thm4.2}, and this manifests $m\left( J_\lambda ,u_\lambda
\right) \geq k$.

We remain to show $\mu _{k+1}\left( A_{\lambda ,u_\lambda }\right) >0$.
Indeed, if the thesis is false, according to the hypothesis $\left(
g_8\right) $,
\begin{eqnarray}
0 &\geq &\mu _{k+1}\left( A_{\lambda ,u_\lambda }\right) =\sup\limits_{%
\stackrel{\frak{B}_k\subset H^2}{\dim \frak{B}_k=k}}\inf\limits_{\stackrel{%
\phi \in \frak{B}_k^{\bot }\cap H^2}{\left\| \phi \right\| _{L^2}=1}%
}\left\langle A_{\lambda ,u_\lambda }\phi ,\phi \right\rangle _{L^2}
\nonumber \\
\ &\geq &\sup\limits_{\stackrel{\frak{B}_k\subset H^2}{\dim \frak{B}_k=k}%
}\inf\limits_{\stackrel{\phi \in \frak{B}_k^{\bot }\cap H^2}{\left\| \phi
\right\| _{L^2}=1}}\left\langle \left( -\Delta +V-\lambda \right) \phi ,\phi
\right\rangle _{L^2}.  \label{eq131}
\end{eqnarray}
As $\mu _{k+1}\left( A_{\lambda ,u_\lambda }\right) <\inf \sigma _{\text{ess}%
}\left( A_{\lambda ,u_\lambda }\right) =-g_0$, once again by Theorem\ref
{thm4.2}, we provoke a ridiculous deduction
\begin{equation}
0\geq \mu _{k+1}\left( A_{\lambda ,u_\lambda }\right) >\mu _{k+1}\left(
A\right) -\sigma _0=0  \label{eq132}
\end{equation}
from $\left( \ref{eq131}\right) $. Thereby we have access to the assertion $%
\left( \ref{eq122}\right) $.

In like manner, we can verify
\begin{equation}
m\left( J_{\sigma _0},u^{*}\right) =M\left( J_{\sigma _0},u^{*}\right) =k.
\label{eq133}
\end{equation}
Let
\[
\Phi \left( w,\tau \right) :=-\Delta w+V\left( x\right) w-\left( \tau
+\sigma _0\right) w-\tau u^{*}-g\left( x,w+u^{*}\right) +g\left(
x,u^{*}\right)
\]
for $\left( w,\tau \right) \in X\times \Bbb{R}$ and clearly $\Phi \in
C^1\left( X\times \Bbb{R},X\right) $. Notice that $\left( w,\tau \right)
=\left( 0,0\right) $ solves $\Phi \left( w,\tau \right) =0$ and $\ker \Phi
_w^{\prime }\left( 0,0\right) =\left\{ 0\right\} $, invoking Theorem2.7.2%
\cite{Niren1}, for suitably small $\delta >0$, $\forall \tau \in \left[
-\delta ,\delta \right] $, there exists a unique $C^1$ solution $w=u\left(
\tau \right) \in X$, s.t.
\begin{equation}
\Phi \left( u\left( \tau \right) ,\tau \right) =0,  \label{eq134}
\end{equation}
where $w=u-u^{*}$, $\tau =\lambda -\sigma _0$. The proof of Theorem\ref
{thm1.1} is complete.\qed\vskip 5pt

\begin{example}
\label{Ex5.1} Take $\alpha >0$ and set
\[
g\left( t\right) =\left\{
\begin{array}{cc}
\alpha \ln \frac 1{1+t}, & t>0, \\
0, & t=0, \\
-\alpha \ln \frac 1{1-t}, & t<0,
\end{array}
\right.
\]
s.t., $\sigma _0-\alpha \leq \mu _i$ for some $i\in \Bbb{N}$, $i\leq k$.
Then the hypotheses $\left( g_1\right) $-$\left( g_8\right) $ holds.
\end{example}

We present a more general assumption replacing $\left( g_4\right) $:

$\left( \widetilde{g}_4\right) $ $g_0<\frac{g\left( x,t\right) }t$ for $x\in
\Bbb{R}^N$, $t\in \Bbb{R\backslash }\left\{ 0\right\} $.

Also suppose:

$\left( \widetilde{g}_6\right) $ $\lambda +g_0>\mu _k$.

Our consequence concerning the nonexistence of nonzero solution reads:

\begin{theorem}
\label{thm5.2} Given $\lambda \in \Bbb{R}$ and $\lambda \leq \sigma _0$,
under the hypotheses $\left( g_2\right) \left( \widetilde{g}_4\right) \left(
\widetilde{g}_6\right) $, $\left( \ref{eq1}\right) $ has no nontrivial
solution.
\end{theorem}

\textit{Proof. }By way of negation, there exists $u_\lambda \neq 0$ solving $%
\left( \ref{eq1}\right) $, i.e., $0\in \sigma _{\text{dis}}\left( \widetilde{%
A}_{\lambda ,u_\lambda }\right) $. Denote
\[
\mu _i\left( \widetilde{A}_{\lambda ,u_\lambda }\right) :=\sup\limits_{%
\stackrel{\frak{B}_{i-1}\subset H^2}{\dim \frak{B}_{i-1}=i-1}}\inf\limits_{%
\stackrel{\phi \in \frak{B}_{i-1}^{\bot }\cap H^2}{\left\| \phi \right\|
_{L^2}=1}}\left\langle \widetilde{A}_{\lambda ,u_\lambda }\phi ,\phi
\right\rangle _{L^2}.
\]
We will show $0\in \rho \left( \widetilde{A}_{\lambda ,u_\lambda }\right) $
and this indeed alludes to the conclusion. Due to $\sigma _{\text{ess}%
}\left( \widetilde{A}_{\lambda ,u_\lambda }\right) =\sigma _{\text{c}}\left(
\widetilde{A}_{\lambda ,u_\lambda }\right) =\left[ \sigma _0-\lambda
-g_0,\infty \right) $, and also notice that $\lambda +\frac{g\left(
x,u_\lambda \right) }{u_\lambda }>\mu _k$ by combining $\left( \widetilde{g}%
_4\right) $ and $\left( \widetilde{g}_6\right) $, if $i\leq k$, exploiting
Theorem\ref{thm4.2} we have
\begin{equation}
\mu _i\left( \widetilde{A}_{\lambda ,u_\lambda }\right) <\sup\limits_{%
\stackrel{\frak{B}_{i-1}\subset H^2}{\dim \frak{B}_{i-1}=i-1}}\inf\limits_{%
\stackrel{\phi \in \frak{B}_{i-1}^{\bot }\cap H^2}{\left\| \phi \right\|
_{L^2}=1}}\left\langle \left( -\Delta +V-\mu _k\right) \phi ,\phi
\right\rangle _{L^2}\leq 0,  \label{eq135}
\end{equation}
and if $i=k+1$, obviously $\mu _i\left( \widetilde{A}_{\lambda ,u_\lambda
}\right) >0$ with the aid of Theorem\ref{thm4.2}. We therefore conclude the
assertion as desired.\qed\vskip 5pt

\textbf{Acknowledgments}. The author would be sincerely indebted to his
family's work behind the scenes over the years, and also appreciate for
strong supports by NSFC(11871066) and E3550105.

\textbf{Conflict of interest Statement}. The author declares that he does
not have any commercial or associative interest that represents a conflict
of interest in connection with the work submitted.

\textbf{Data Availability Statement.} The manuscript has no associate data.

\section{Appendix}

\begin{theorem}
\label{thm6.1} Under the hypotheses of Theorem\ref{thm2.2}, suppose $\frac
N2<p\leq N$ with $N\geq 4$, or $N\leq 3$. Let $\lambda _n\rightarrow \lambda
_0$ and $\left\{ u_{\lambda _n}\right\} _{n=1}^\infty $ be the solution
sequence of $\left( \ref{eq1}\right) $ for $\lambda =\lambda _n$, s.t. $%
\left\{ u_{\lambda _n}\right\} _{n=1}^\infty $ is bounded in $X$, then $%
\left\{ u_{\lambda _n}\right\} _{n=1}^\infty $ does as well in $L^\infty
\left( \Bbb{R}^N\right) $.
\end{theorem}

\textit{Proof. }According to the hypothesis, $\exists M>0$, $\left\|
u_{\lambda _n}\right\| _{H^1}\leq M$. Invoking Theorem\ref{thm2.2}, $\left\{
u_{\lambda _n}\right\} _{n=1}^\infty \subset H^2\left( \Bbb{R}^N\right) \cap
L^\infty \left( \Bbb{R}^N\right) $. Notice that
\begin{eqnarray}
\int_{\Bbb{R}^N}\left| \Delta u_{\lambda _n}\right| ^2dx &\leq &\int_{\Bbb{R}%
^N}\left( \left| V_1\right| +\left\| V_2\right\| _{L^\infty }+\left| \lambda
_n\right| +C\right) ^2\cdot \left| u_{\lambda _n}\right| ^2dx  \nonumber \\
&\leq &2\int_{\Bbb{R}^N}\left| V_1\right| ^2\left| u_{\lambda _n}\right|
^2dx+2\left( \left\| V_2\right\| _{L^\infty }+\left| \lambda _n\right|
+C\right) ^2\int_{\Bbb{R}^N}\left| u_{\lambda _n}\right| ^2dx  \nonumber \\
&&\   \label{eq136}
\end{eqnarray}
and
\begin{eqnarray}
\int_{\Bbb{R}^N}\left| V_1\right| ^2\left| u_{\lambda _n}\right| ^2dx &\leq
&\left[ \int_{\Bbb{R}^N}\left| V_1\right| ^{2\widetilde{p}}dx\right] ^{\frac
2{\widetilde{p}}}\cdot \left[ \int_{\Bbb{R}^N}\left| u_{\lambda _n}\right|
^{2\widetilde{q}}dx\right] ^{\frac 1{\widetilde{q}}}  \nonumber \\
\ &=&\left\| V_1\right\| _{L^p}^4\cdot \left[ \int_{\Bbb{R}^N}\left|
u_{\lambda _n}\right| ^{2\widetilde{q}}dx\right] ^{\frac 1{\widetilde{q}}}
\nonumber \\
\ &\leq &\left\| V_1\right\| _{L^p}^4\cdot \left\| u_{\lambda _n}\right\|
_{L^\infty }^{2-\frac{2^{*}}{\widetilde{q}}}\cdot \left\| u_{\lambda
_n}\right\| _{L^{2^{*}}}^{\frac{2^{*}}{\widetilde{q}}}  \nonumber \\
\ &\leq &\left\| V_1\right\| _{L^p}^4\cdot \left\| u_{\lambda _n}\right\|
_{L^\infty }^{2-\frac{2^{*}}{\widetilde{q}}}\cdot C\left( 1,2,N\right) ^{%
\frac{2^{*}}{\widetilde{q}}}\cdot \left\| u_{\lambda _n}\right\| _{H^1}^{%
\frac{2^{*}}{\widetilde{q}}}  \nonumber \\
\ &\leq &M^{\frac{2^{*}}{\widetilde{q}}}\left\| V_1\right\| _{L^p}^4C\left(
1,2,N\right) ^{\frac{2^{*}}{\widetilde{q}}}\left\| u_{\lambda _n}\right\|
_{L^\infty }^{2-\frac{2^{*}}{\widetilde{q}}}.  \label{eq137}
\end{eqnarray}
$\widetilde{p}=\frac p2$, $\frac 1{\widetilde{p}}+\frac 1{\widetilde{q}}=1$,
$2^{*}<2\widetilde{q}<\frac{2N}{N-4}$. Therefore
\begin{eqnarray}
\int_{\Bbb{R}^N}\left| \Delta u_{\lambda _n}\right| ^2dx &\leq &2M^{\frac{%
2^{*}}{\widetilde{q}}}\left\| V_1\right\| _{L^p}^4C\left( 1,2,N\right) ^{%
\frac{2^{*}}{\widetilde{q}}}\left\| u_{\lambda _n}\right\| _{L^\infty }^{2-%
\frac{2^{*}}{\widetilde{q}}}+2\left( \left\| V_2\right\| _{L^\infty }+\left|
\lambda _n\right| +C\right) ^2M^2.  \nonumber \\
&&  \label{eq138}
\end{eqnarray}
Take $\varepsilon \leq \frac 1{2K\left( 2,2,N\right) }$ and thus
\begin{eqnarray}
\left\| \nabla u_{\lambda _n}\right\| _{L^2} &\leq &\left\| \Delta
u_{\lambda _n}\right\| _{L^2}+\left[ 1+2K\left( 2,2,N\right) \varepsilon
^{-1}\right] \cdot \left\| u_{\lambda _n}\right\| _{L^2}  \nonumber \\
\ &\leq &2^{\frac 12}M^{\frac{2^{*}}{2\widetilde{q}}}\left\| V_1\right\|
_{L^p}^2C\left( 1,2,N\right) ^{\frac{2^{*}}{2\widetilde{q}}}\left\|
u_{\lambda _n}\right\| _{L^\infty }^{1-\frac{2^{*}}{2\widetilde{q}}%
}+2^{\frac 12}\left( \left\| V_2\right\| _{L^\infty }+\left| \lambda
_n\right| +C\right) M  \nonumber \\
&&\ \ +\left[ 1+2K\left( 2,2,N\right) \varepsilon ^{-1}\right] M.
\label{eq139}
\end{eqnarray}
This yields
\begin{eqnarray}
\left\| u_{\lambda _n}\right\| _{H^2} &\leq &2^{\frac 32}M^{\frac{2^{*}}{2%
\widetilde{q}}}\left\| V_1\right\| _{L^p}^2C\left( 1,2,N\right) ^{\frac{2^{*}%
}{2\widetilde{q}}}\left\| u_{\lambda _n}\right\| _{L^\infty }^{1-\frac{2^{*}%
}{2\widetilde{q}}}  \nonumber \\
&&\ +\left[ 2^{\frac 32}\left( \left\| V_2\right\| _{L^\infty }+\left|
\lambda _n\right| +C\right) +2K\left( 2,2,N\right) \varepsilon
^{-1}+2\right] M.  \label{eq140}
\end{eqnarray}
Combining $\left( \ref{eq54}\right) $,$\left( \ref{eq55}\right) $,$\left(
\ref{eq56}\right) $ and $\left( \ref{eq57}\right) $, there exists a $C^{*}>0$%
, s.t.
\begin{equation}
\left\| u_{\lambda _n}\right\| _{L^\infty }\leq C^{*}\left\| u_{\lambda
_n}\right\| _{H^2}.  \label{eq141}
\end{equation}
Inserting $\left( \ref{eq140}\right) $ into $\left( \ref{eq141}\right) $, we
conclude the proof for the case $N\geq 5$. Analogously we can verify the
case $N\leq 4$.\qed


\vskip 5pt

\begin{thebibliography}{99}
\bibitem{AF}  Adams, R.A., Fournier, J.J.F.: Sobolev Spaces, volume 140 of
Pure and Applied Mathematics. Elsevier-Academic Press, Amsterdam(2003)

\bibitem{ADN1}  Agmon, S., Douglis, A., Nirenberg, L.: Estimates near the
boundary for solutions of elliptic partial differential equations satisfying
general boundary conditions II. Comm.Pure Appl.Math. 17, 35-92(1964)

\bibitem{ADN2}  Agmon, S., Douglis, A., Nirenberg, L.: Estimates near the
boundary for solutions of elliptic partial differential equations satisfying
general boundary conditions I. Comm.Pure Appl.Math. 12, 623-727(1959)

\bibitem{B}  Brezis, H.: Analyse Fonctionelle, Th\'eorie et Applications.
Masson, Paris(1983)

\bibitem{Cha}  Chang, K.C.: Infinite dimensional Morse theory and multiple
solution problem. Birkh\"auser, Boston(1993)

\bibitem{CT}  Costa, D., Tehrani, H.: On a class of asymptotically linear
elliptic problems in $\Bbb{R}^N$. J.Differential Equations. 173,
470-494(2001)

\bibitem{Heinz}  Heinz, H.P.: On the number of solutions of nonlinear
Schr\"odinger equations and on unique continuation. J.Differential
Equations. 116, 149-171(1995)

\bibitem{HS}  Hislop, P.D., Sigal, I.M.: Introduction to spectral theory
with application to Schr\"odinger operators. Springer, New York(1996)

\bibitem{J}  Jeanjean, L.: On the existence of bounded Palais-Smale
sequences and applications to a Landesman-Lazer-type problem set on $\Bbb{R}%
^N$. Proc.Roy.Soc.Edinburgh Sect.A. 129, 787-809(1999)

\bibitem{JT}  Jeanjean, L., Tanaka, K.: A positive solution for an
asymptotically linear elliptic problem on $\Bbb{R}^N$ autonomous at
infinity. ESAIM Control.Optim.Calc.Var. 7, 597-614(2002)

\bibitem{JK}  Jerison, D., Kenig, C.E.: Unique continuation and absence of
positive eigenvalues for Schr\"odinger operators. Ann.of Math. 121,
463-494(1985)

\bibitem{K}  Kenig, C.E.: Carleman estimates, uniform Sobolev inequalities
for second-order differential operators, and unique continuation theorems
ICM Proceedings. Berkeley, California(1986)

\bibitem{LiLi}  Li, C., Li.S.J.: The Fu\v c\'\i k spectrum of Schr\"odinger
operator and the existence of four solutions of Schr\"odinger equations with
jumping nonlinearities, J.Differential Equations. 263, 7000-7097(2017)

\bibitem{LZ}  Li, G.B. Zhou, H.S.: The existence of a positive solution to
asymptotically linear scalar field equations. Proc.Roy.Soc.Edinburgh Sect.A.
130, 81-105(2000)

\bibitem{Lions}  Lions, P.L.: The concentration-compactness principle in the
calculus of variations. The locally compact case. Ann.Inst.Henri Poincar\'e,
Analyse Non Lin\'eaire. 1, 109-145, 223-283(1984)

\bibitem{MS}  Maia, L.A., Soares, M.: Spectral theory approach for a class
of radial indefinite variational problems. J. Differential
Equations.(2018).https://doi.org/10.1016/j.jde.2018.11.020

\bibitem{Niren1}  Nirenberg, L.: Topics in Nonlinear Functional Analysis.
Courant Lect. Notes Math. New York University(1974)

\bibitem{Niren2}  Nirenberg, L.: Uniqueness in Cauchy problems for
differential equations with constant leading coefficients. Comm.Pure
Appl.Math.. 10, 89-105(1957)

\bibitem{RS}  Reed, M., Simon, B.: Methods of Modern Mathematical Physics,
vol.4, Analysis of Operators. Academic Press, New York(1978)

\bibitem{S}  Song, L.J.: Existence andmultiplicity of solutions for
Schr\"odinger equations with asymptotically linear nonlinearities allowing
interaction with essential spectrum. Partial Differential Equations and
Applications.(2022).https://doi.org/10.1007/s42985-022-00162-7

\bibitem{SZ}  Stuart, C.A., Zhou, H.S.: Applying the mountain pass theorem
to an asymptotically linear elliptic equation on $\Bbb{R}^N$, Comm.Partial
Differential Equations. 24, 1731-1758(1999)

\bibitem{W}  Watanabe, T.: Radial solutions with a vortex to an
asymptotically linear elliptic equation. NoDEA Nonlinear Differential
Equations Appl. 15, 387-411(2008)

\bibitem{Willem}  Willem, M.: Minimax theorems, Progress in Nonlinear
Differential Equations and Their Applications, 24. Birkh\"auser, Boston(1996)
\end{thebibliography}
\end{document}